\documentclass[11pt,a4paper]{article}
\usepackage[english]{babel}
\usepackage[latin1]{inputenc}
\usepackage{amssymb}
\usepackage{amsmath}
\usepackage{amsthm}
\usepackage{bbm}

\newcommand\RR{\mathbb{R}}
\newcommand\RRB{\mathbb{R} \times B}
\newcommand\tr{\ensuremath{\triangle}}
\newcommand\Str{Strichartz }
\newcommand\Schr{Schrödinger }
\newcommand\jnk{\ensuremath{j_{n+\frac 12,k}}}
\newcommand\Ps{\ensuremath{\mathcal{P}_s}}
\newcommand\Psz{\ensuremath{\mathcal{P}_{s_0}}}
\newcommand\Gs{\ensuremath{\mathcal{G}_s}}
\newcommand\Gsz{\ensuremath{\mathcal{G}_{s_0}}}

\newcommand\norm[2]{\ensuremath{|\!| #1 |\!|_{#2}}}

\newcommand\nLi[3]{\norm{#1}{L^#2(#3)}}
\newcommand\nLdi[2]{\norm{#1}{L^2(#2)}}

\newcommand\nH[2]{\norm{#1}{H^{#2}}}

\newcommand\nX[2]{\norm{#1}{#2}}

\newcommand\nXsb[3]{\norm{#1}{X^{#2,#3}(\RR \times B)}}
\newcommand\nXsbT[3]{\norm{#1}{X_T^{#2,#3}(B)}}

\newtheorem{Def}{Definition}
\newtheorem{thm}{Theorem}[section]
\newtheorem{lem}[thm]{Lemma}
\newtheorem{prop}[thm]{Proposition}
\newtheorem{rem}{Remark}
\newtheorem{cor}[thm]{Corollary}

\newtheorem*{citelem}{Lemma}

\title{Cubic nonlinear Schrödinger equation on three dimensional balls with radial data}
\author{Ramona Anton}
\date{}
\begin{document}

\maketitle

\begin{abstract}

We prove wellposedness of the Cauchy problem for the cubic nonlinear Schrödinger equation with Dirichlet boundary conditions and radial data on 3D balls. The main argument is based on a bilinear eigenfunction estimate and the use of $X^{s,b}$ spaces. The last part presents a first attempt to study the non radial case. We prove bilinear estimates for the linear Schrödinger flow with particular initial data.

\end{abstract}

\section{Introduction}
Let us denote by $B$ the unit ball in $\RR^3$. We are interested in the Dirichlet problem for the semilinear Schrödinger equation on $B$. We denote by $\tr = \tr_D$ the Dirichlet Laplacian on $B$.
\begin{equation}
\label{NLS}
\left \{
\begin{array}{rcl}
i \partial_t u + \triangle u &=& |u|^2 u, \ on \ \RRB \\
u_{|_{t=0}} &=& u_0, \ on \ B\\
u_{|_{\mathbb R \times \mathbb{S}^2}} &=& 0.
\end{array}
\right.
\end{equation}

In the Euclidean context the nonlinear Schrödinger equation has been extensively studied since the seventies. One of the main tools in studying local existence is the generalized Strichartz inequality. We could mention the work of \Str \cite{STR}, Ginibre-Velo \cite{GiVe}, Keel-Tao \cite{KT}, Cazenave-Weissler \cite{CaW}, Yajima \cite{Ya} and Tsutsumi \cite{Tsu}.
In the last decade there has been an intense activity on studying the influence of the geometry on the behavior of Schrödinger flow. Recent results have shown that the geometry plays a major role : see for example the wellposedness result in $H^s$ on a square of $  \RR^2$, Bourgain \cite{BourgBook}, for $s>0$, as opposed to an illposedness result in $H^s$ on a disc of $\RR^2$, Burq-Gérard-Tzvetkov \cite{BGTgafa}, for $s<\frac 13$. For cubic equations a general results on domains of $\RR^2$ is due to Brezis-Gallouet \cite{BrGa} and Vladimirov \cite{Vl}. This result is based on energy methods and logarithmic inequalities. It does not provide informations on the regularity of the flow.

Another direction is to prove Strichartz inequality with loss of derivatives. For example on boundaryless compact manifolds or asymptotically flat metrics Strichartz inequality with loss of derivatives were proved by Staffilani-Tataru \cite{StTa} and Burq-Gérard-Tzvetkov \cite{BGTajm}. Moreover, it was shown that in some geometries a loss of derivatives is inevitable. 

For the boundaryless case the bilinear Strichartz estimates proved to be a useful tool.  They were used by several authors in the context of wave and \Schr equations : Klainerman-Machedon \cite{KlMa}, Klainerman-Machedon-Bourgain-Tataru \cite{KlMaBourgTa}, Foschi-Klainerman \cite{FoKl}, Bourgain \cite{Bourg1,Bourg2}, Tao \cite{Tao}, Burq-Gérard-Tzvetkov \cite{BGTens, BGTinv} and references therein.

As we show in the sequel, we used them successfully also in the case of radial data on balls of $\RR^3$. Notice that in the case of a boundary domain the gradient does not preserve the intersection of domains of $(-\tr_D)^k$, $k\in \mathbb{N}$, as opposed to the boundaryless case. Therefore, we introduce a supplementary bilinear estimates, in which we have a gradient on one of the terms. This property helps us to handle integration by parts and appears to be a useful complement to the usual bilinear estimates in boundary domains. Bilinear estimates have the advantage of showing interactions between the large and the small frequencies, which is extremely useful when treating the non linear terms.

\begin{rem}
\label{deltaNo}
For technical reasons, in Section \ref{PsZn1}, which treats the case of a special nonradial initial data, we need to use numbers $N\in \delta^{\mathbb{N}}$, where $1<\delta$ and close to $1$ (for example $\delta<\sqrt{1.5}$) instead of the usual powers of $2$. For the coherence of the paper we shall consider hereafter $N\in \delta^{\mathbb{N}}$. We introduce the notation $n \sim_\delta N$ to denote $n\in [N,\delta N]$.
\end{rem}

\begin{Def} Let $s>0$. We say that $S(t)=e^{it\tr_D}$, the flow of the linear \Schr equation with Dirichlet boundary conditions on the ball $B\subset \RR^3$, satisfies property $\Ps$ if for all $N,L>0$ and $u_0, v_0 \in L^2(B)$ satisfying 
\begin{equation}
\label{scoNL}
\mathbbm{1}_{\sqrt{-\tr} \in {[N,\delta N]}} u_0=u_0,\ \ \mathbbm{1}_{\sqrt{-\tr} \in {[L,\delta L]}} v_0=v_0,
\end{equation}
the following estimate holds :
\begin{equation}
\label{bilinEst}
	\nLdi{S(t)u_0 S(t)v_0}{(0,1)\times B} \leq C (\min(N,L))^{s} \nLdi{u_0}{B}\nLdi{v_0}{B}.
\end{equation}

\noindent We say that $S(t)$ satisfies property $\Gs$ if for $u_0 , v_0 \in H_0^1(B)$ satisfying (\ref{scoNL}) : \begin{equation}
\label{dISB}
	\nLdi{(\triangledown S(t)u_0) S(t)v_0}{(0,1)\times B} \leq C N(\min(N,L))^{s} \nLdi{u_0}{B}\nLdi{v_0}{B}.
\end{equation}
\end{Def}

For $s\geq 0$ we denote by $H_D^s(B)$ the domain of $(-\tr_D)^\frac s2$.  Let us remind that for $s\in [0,\frac 12)$ we have $H_D^s(B)=H^s(B)$, for $s\in (\frac 12, 1]$ we have $H_D^s(B)=H_0^s(B)$ and for $1\leq s \leq 2$, $H_D^s(B)=H^s(B) \cap H_0^1(B)$.

\begin{thm} \label{lThm}Let $s_0>0$ and assume that $S(t)$ satisfies properties $\Psz$ and $\Gsz$. Then for $s>s_0$ and for $R>0$ there exists $T=T(R,s)>0$ and $X(R,s) \subset C([-T,T], H_D^{s}(B))$ such that for every $u_0\in H_D^{s}(B)$ with $\nX{u_0}{H_D^{s}(B)} \leq R$ there exists a unique solution $u\in X(R,s)$ of the Cauchy problem (\ref{NLS}). If moreover $u_0\in H_D^{\sigma}(B)$, for some $\sigma>s_0$, then the regularity propagates, i.e. $u \in C([-T,T], H_D^{\sigma}(B)$. The flow $u_0\in B_{H_D^{s}(B)}(0,R) \mapsto u\in C([-T,T],H_D^{s}(B))$ is Lipschitz on bounded subsets of $H_D^{s}(B)$.
\end{thm}
The proof of Theorem \ref{lThm} relies on the classical use of $X^{s,b}$ spaces and a contraction mapping argument. However, the new element is estimate $\Gs$ which we use to treat the case when one of the frequencies is much larger in front of the others. Having a general theorem, it suffices to prove properties $\Ps$ and $\Gs$ in order to obtain a local existence result. Notice that, for $0<s_0<1$, we obtain a local existence result in $H_0^1(B)$, which combined with the conservation of the $H_0^1(B)$ norm (energy conservation) implies global existence. However, at this moment we are only able to prove properties $\Ps$ and $\Gs$, with $s<1$, in the radial case. In the case of radial data, we define $$L_{rad}^2(B) = \{f\in L^2(B)\ s.t.\ \exists g:\RR\rightarrow \mathbb{C} ,\ f(x)=g(|x|)\}.$$
\begin{prop} \label{bilinLem}
Let $s >\frac 12$. For $u_0,v_0 \in L_{rad}^2$, $S(t)$ satisfies (\ref{bilinEst}) (i.e. satisfies $\Ps$ for radial data) and  (\ref{dISB}) (i.e. $\Gs$ for radial data).
\end{prop}

Notice that the radial symmetry is conserved by the non-linear flow, therefore having properties $\Ps$ and $\Gs$ for radial data in Lemma \ref{bilinLem} implies the wellposedness for radial data in $H_D^s(B)$. This fact combined with the conservation of energy gives a global existence result for the defocusing cubic NLS on $B$.

\begin{cor} For $u_0 \in L_{rad}^2(B)\cap H_0^{1}(B)$ there exists a global solution $u\in C((-\infty, \infty), H_0^1(B))$ to (\ref{NLS}), satisfying uniqueness in $X^{1,b}(B)$ for some $b>\frac 12$, regularity property : if $u_0\in H^2(B)\cap H_0^1(B)$ then $u\in C(\RR,H^2(B)\cap H_0^1(B))$ and $C^\infty$ and Lipschitz regularity on bounded subsets of $H_D^s(B)$ of the flow $u_0\in u$.
\end{cor}

Notice that a loss of $\frac{1}{2}$ derivative appears in the end point Strichartz estimate of the linear \Schr flow with radial data on $B$ (for definitions see for example \cite{KT}). We test it against an eigenfunction of $\tr$ on $B$ : $e_n(r)=\frac{\sin(n\pi r)}{r}$. 
Estimating $e^{it\tr} e_n(r) = e^{i t n^2}\frac{\sin(n\pi r)}{r} $ in the  $L^2_t(I, L^6_x(B))$ norm gives $$\nLi{e^{it\tr} e_n}{2}{I,L^6(B)} \leq c \nLi{\frac{\sin(n\pi r)}{r}} {6}{B} = c_n n^{\frac 12} = c c_n \nX{e_n}{H^\frac 12(B)},$$ where $c_n\rightarrow c_0$ when $n\rightarrow \infty$.

Let us mention that for general data on the ball $B$, property $\Ps$ is false for $s\leq \frac{7}{12}$ (see \cite{BGTcortona}). Indeed, an counter example using an eigenfunction associated to the first zero of the Bessel function has been constructed in that paper.

One step in proving $\Ps$ for $S(t)$ acting on general data is to prove it for initial data in the vector space spanned by the eigenfunctions corresponding to the first zero of Bessel function.

\begin{prop}
\label{ISBzn1}
The linear flow $S(t)=e^{it\tr}$ satisfies property $\Ps$ for $s>\frac {11}{12}$ and for initial data in the vector space spanned by the eigenfunctions corresponding to the first zero of Bessel function. With the notations of Section \ref{PsZn1}, this reads
$$\nLdi{S(t)u_{0N}S(t)v_{0N}}{(0,1)\times B} \leq c min(N,L)^{\frac {11}{12}+\epsilon} \nLdi{u_{0N}}{B} \nLdi{v_{0N}}{B}.$$
\end{prop}

The proof of Proposition (\ref{ISBzn1}) uses some tools from analytic number theory and precised asymptotics of the first positif zero of the Bessel function. We are able to prove that $\Ps$ holds on such data for an $s<1$, which is encouraging in the perspective of proving the existence of global strong solutions for (\ref{NLS}). 
However, such data are not conserved by the nonlinear flow and therefore at this stage we cannot apply Theorem \ref{lThm}.

This paper is organized as follows : in Section \ref{Xsb} we introduce $X^{s,b}(B)$ spaces and present the proof of Theorem \ref{lThm}. In Section \ref{bilinear} we give the proof of properties $\Ps$ and $\Gs$ in the radial case, as well as the proof of Proposition \ref{ISBzn1}. In Appendix we present the proofs of some technical lemmas. 

\noindent {\bf Acknowledgments :} \textit{The author would like to thank E.Fouvry for valuable discussions on exponential sums. She is also indebted to P.Gérard for suggesting the problem and helping to its achievement. This result is part of author's PhD thesis in preparation at Université Paris Sud, Orsay, under P.Gérard's direction.}

\section{Proof of Theorem \ref{lThm}}
\label{Xsb}
In this section we introduce the $X^{s,b}$ Bourgain type spaces on $B$. We present the classical (by now) way of proving local wellposedness for the cubic equation from bilinear estimates (\ref{bilinEst}). We shall prove those estimates in Section \ref{bilinear}.

\subsection{$X^{s,b}$ spaces}

They are spaces of functions in time and space variables, introduced by Bourgain for the Schrödinger operator. On $\RR^d$ the definition is given in terms of multipliers of the Fourier transform and $L^2$ spaces. We shall follow  the definition given by Burq, Gérard and Tzvetkov \cite{BGTinv} using spectral projectors on manifolds.

Using the notation $\langle x \rangle=(1+x^2)^{\frac 12}$, we have the following definition of $X^{s,b}$ spaces :

\begin{Def} Let $(e_n)_n$ be a $L^2$ orthonormal basis of eigenfunctions of the Dirichelet Laplacian $-\tr$ with eigenvalues $\mu_n$ on $B$. Let $\Pi_n$ be the orthogonal projector along $e_n$. Then, for $s\geq 0$ and $b\in \RR$, $$X^{s,b}(\RR\times B)= \{u\in \mathcal{S}'(\RR,L^2(B))\ s.t. \ \nXsb{u}{s}{b} <\infty\},$$ 
where we denote by $\nXsb {u}{s}{b}$ the norm 
\begin{equation}
\label{XsbDef}
\nXsb {u}{s}{b} ^2=\sum_n \nLi{ \langle\tau + \mu_n\rangle^{\frac b2} \langle\mu_n\rangle^{\frac s2} \widehat{\Pi_n u}(\tau)}{2}{\RR_\tau, L^2(B)}^2,
\end{equation}
and $\widehat{\Pi_n u}(\tau)$ denotes the Fourier transform of $\Pi_n u$ with respect to the time variable. Moreover, for $u\in X^{0,\infty}(B) = \cap_{b\in \RR} X^{0,b}(B)$, we define, for $s\leq 0$ and $b\in \RR$, the norm $\nXsb {u}{s}{b}$ by (\ref{XsbDef}).
\end{Def}

We use a contraction mapping argument to obtain local existence. We therefore need to define some local in time version of $X^{s,b}(\RRB)$. For $T>0$ we denote by $X_T^{s,b}(B)$ the space of restrictions of elements of $X^{s,b}(\RR \times B)$ endowed with the norm
$$\nX{u}{X_T^{s,b}} = \inf \{\nX{\tilde u}{X^{s,b}(\RR,B)},\ \tilde u_{|_{]-T,T[\times B}}=u \}.$$

Using the notation $S(t)=e^{it\tr}$ for the linear \Schr flow, we prove the following property, which will we used in order to estimate the $X^{s,b}$ norm.
\begin{prop} Let $s\geq 0$ and $u\in \mathcal{S}'(\RR, L^2(B))$. 
Then we have the following equivalence $u\in X^{s,b}(\RR \times B)$ $\iff$ $S(-t)u(t,\cdot) \in H^b(\RR, H_D^s(B))$ and 
\begin{equation}
\label{equivDef}
\nXsb{u}{s}{b} = \nX{S(-t)u(t,\cdot)}{H^b(\RR,H^s(B))}.
\end{equation}

Moreover, for $b>\frac 12$, $X^{s,b}(\RR\times B) \hookrightarrow C(\RR, H_D^s(B))$.
\end{prop}

\begin{proof} Let us denote by $F(t,\cdot)=S(-t)u(t,\cdot)$. Then $F\in \cal{S}'(\RR \times B)$ and $\Pi_n F (t,\cdot)= e^{it\mu_n} \Pi_n u(t,\cdot)$. Consequently, $\widehat{\Pi_n F}(\tau) = \widehat{\Pi_n u} (\tau-\mu_n)$. Introducing this identity in (\ref{XsbDef}), we conclude that $\nXsb{u}{s}{b} = \nX{F}{H^b(\RR, H^s(B))}$.

In the case $b>\frac 12$,  $H^b(\RR, H_D^s(B))\hookrightarrow C(\RR, H_D^s(B))$ and since $u(t,\cdot)=S(t)F(t,\cdot)$, we deduce $u\in C(\RR, H^s(B))$.
\end{proof}

Since for $b>\frac 12$ the embedding $X^{s,b}(B) \hookrightarrow C(\RR, H_D^s(B))$ holds and continuity is a local property, we deduce that $X_T^{s,b}(B) \hookrightarrow C((-T,T), H_D^s(B))$. Therefore, we shall focus on proving a local existence theorem in a $X_T^{s,b}(B)$ space for some $b>\frac 12$.

\begin{rem} From the Sobolev embedding $H^{\frac 14}(\RR) \hookrightarrow L^4(\RR)$ we deduce :
\begin{equation}
\label{X014}
X^{0, \frac 14}(\RR\times B) \hookrightarrow L^4(\RR, L^2(B)).
\end{equation}
\end{rem}
\noindent Indeed, using the conservation of the $L^2$ norm by the linear \Schr flow $e^{it\tr}$ and the definition (\ref{equivDef}) of the $X^{s,b}$ norm,
$$\nLi{f}{4}{\RR,L^2(B)} = \nLi{e^{it\tr}f}{4}{\RR,L^2(B)} \leq \nX{e^{it\tr}f}{H^{\frac 14}(\RR, L^2(B))} = \nXsb{f}{0}{\frac 14}.$$

Let us give the bilinear estimates (\ref{bilinEst}) and (\ref{dISB}) in the $X^{s,b}$ context. 

\begin{lem} 
Let $s>0$ and $N,L >0$. 
If $S(t)$ satisfies property $\Ps$ then for every $b>\frac 12$ there exists $c_b>0$ such that, for $f,g \in L^2(B)$ satisfying (\ref{scoNL}),
\begin{equation}
	\label{f2ISB}
	\nLdi{f g}{\RR \times B} \leq c_b (\min(N,L))^s \nXsb{f}{0}{b}\nXsb{g}{0}{b}.
\end{equation}

If $S(t)$ satisfies property $\Gs$ then for $b>\frac 12$ there exists $c_b>0$ such that for $f,g \in X^{0,b}(\RR \times B)$ satisfying (\ref{scoNL}),
\begin{equation}
	\label{f3ISB}
	\nLdi{(\nabla f) g}{\RR \times B} \leq c_b N (\min(N,L))^s \nXsb{f}{0}{b}\nXsb{g}{0}{b}.
\end{equation}
\end{lem}

\begin{proof}
For the proof of the first property of this "transfer lemma" we shall refer to Lemma 2.3 of \cite{BGTinv}. 

For the second property we follow closely the proof in Lemma 2.3 of \cite{BGTinv}. We suppose first that $f,g\in X^{0,b}(\RRB)$ are supported in $(0,1)\times B$. We write $f(t,\cdot)=S(t)F(t,\cdot)$, where $F(t,\cdot)=S(-t)f(t,\cdot)$ and similarly for $g$. By (\ref{equivDef}), $\nXsb{f}{0}{b} = \nX{F}{H^b(\RR, L^2(B))}$. Using the inverse Fourier transform we write $F(t)=\frac{1}{2\pi} \int_{\RR} e^{it\tau} \hat{F}(\tau){\rm d}\tau$, where $\hat{F}(\tau)$ designs the Fourier transform of $F$ in $t$ variable. Therefore $f(t)=\frac{1}{2\pi} \int_{\RR} e^{it\tau} S(t)\hat{F}(\tau){\rm d}\tau$. Consequently,
$$(\nabla f) g = \frac{1}{(2\pi)^2} \int_{\RR}\int_{\RR} e^{it(\tau + \eta)} \nabla S(t) \hat{F}(\tau) S(t)\hat{G}(\eta){\rm d}\eta {\rm d}\tau.$$ By property $\Gs$, we estimate $$\nLdi{(\nabla f) g}{(0,1)\times B} \leq c N (\min(N,L))^s \int_{\RR}\int_{\RR} \nLdi{\hat{F}(\tau) }{B} \nLdi{\hat{G}(\eta)}{B} {\rm d}\eta {\rm d}\tau.$$ We use the condition $b>\frac 12$ to insure integrability in $\tau$ and $\eta$. By Cauchy Schwarz, the double integral is bounded by $$c \nLdi{\langle \tau \rangle ^b \hat{F}(\tau)}{\RR_\tau \times B} \nLdi{\langle \eta \rangle ^b \hat{G}(\eta)}{\RR_\eta \times B}.$$ Using (\ref{equivDef}) we deduce (\ref{f3ISB}) for $f,g\in X^{0,b}(\RRB)$ supported in $(0,1)\times B$. For the general case we decompose $f(t)=\sum_{k\in \mathbb{Z}} \psi(t-\frac{k}{2}) f(t)$ for $\psi \in C_0^\infty(\RR)$ such that $\sum_{k\in \mathbb{Z}} \psi(t-\frac{k}{2}) = 1$ for all $t\in \RR$. We do the same decomposition for $g$. Using the almost disjoint supports, the general case follows from the particular case of $f,g$ supported in $(0,1)\times B$.
\end{proof}

Those spaces allow to estimate easily the linear flow.

\begin{prop} Let $b,s>0$ and let $u_0 \in H^s(B)$. Then
\begin{equation}
\label{linEst}
\nX{S(t)u_0}{X_T^{s,b}} \leq c T^{\frac 12 -b} \nH{u_0}{s}
\end{equation}
\end{prop}
\begin{proof} Indeed, let $\epsilon>0$ and $\varphi \in C_0^\infty(\RR)$, $\varphi\equiv 1$ on $]-T-\epsilon,T+\epsilon[$. Then 
$$\nXsbT{S(t)u_0}{s}{b} \leq \nXsb{\varphi(t)S(t)u_0}{s}{b} \leq \nX{\varphi(t)u_0}{H^b(\RR,H_D^s(B))} \leq c T^{\frac 12 -b}\nH{u_0}{s}.$$
\end{proof}

Consequently, the difficulty concentrates on proving estimates adapted to the nonlinearity. Using the Duhamel formula for a nonlinear Schrödinger equation, we know that the solution reads, at least formally, $u(t)=S(t)u_0 -i \int_0^t S(t-t') f(t') {\rm d}t'.$ Thus, we have to estimate the $X_T^{s,b}(B)$ norm of 
\begin{equation}
\label{nlTerm}
w(t)=\int_0^t S(t-t') f(t') {\rm d}t'.
\end{equation}
This is treated by the following general lemma due to Bourgain \cite{Bourg2}, but we refer to Ginibre \cite{Gi} for a simpler proof.

\begin{citelem} {\upshape \textbf{(\cite{Bourg2}, \cite{Gi})}} Let $0<b'<\frac 12$ and $0<b<1-b'$. 
Then for all $f\in X_T^{s,-b'}(B)$, having defined $w$ as in (\ref{nlTerm}), we have $w\in X_T^{s,b}(B)$ and moreover
\begin{equation}
\label{nlTermEst}
\nX{w}{X_T^{s,b}(B)} \leq CT^{1-b-b'} \nX{f}{X_T^{s,-b'}(B)}.
\end{equation}
\end{citelem}

Please note that $X^{s,b} \subset X^{s,-b'}$. We shall use the previous lemma to apply a fix point method, as we start from a bigger space $X_T^{s,-b'}$ to arrive in a smaller one, $X_T^{s,b}$.

The nonlinearity we are interested in is the cubic one. Therefore we prove the following fundamental lemma.

\begin{lem} \label{Flem}
Assuming properties $\Ps$ and $\Gs$, for $s>\frac 12 $, there exists $(b,b')\in \RR^2$, satisfying
\begin{equation} 
\label{star}
0<b'<\frac 12 <b,\ b+b'<1,
\end{equation}
and $C>0$ such that for every triple $(u_j)$, $j\in \{1,2,3\}$, in $X^{s,b}(\RR\times B)$,
\begin{equation}
\label{FNLest}
\nXsb{u_1 u_2 u_3}{s}{-b'} \leq c \prod_{j=1}^3 \nXsb{u_j}{s}{b}.
\end{equation}
\end{lem}

Let us show how this lemma implies the local wellposedness Theorem \ref{lThm}.

\begin{proof} (of Theorem \ref{lThm})
Solving the NLS equation (\ref{NLS}) is equivalent with solving the Duhamel integral equation, with Dirichlet boundary conditions :
\begin{equation}
\label{2star}
u(t) = S(t)u_0 -i \int_0^t S(t-\tau) \{|u(\tau)|^2 u(\tau)\}{\rm d}\tau.
\end{equation}
Let us denote by $\Phi(u)$ the left hand side of the equation. 
We consider $(b,b')\in \RR^2$ given by Lemma \ref{Flem}.

Let $R>0$ and $u_0\in H_D^s(B)$ such that $\nH{u_0}{s}\leq R$. We show that there exists $R'>0$ and $0<T<1$ depending on $R$ such that $\Phi$ is a contracting map from the ball $B(0,R') \subset X_T^{s,b}(B)$ onto itself.

From the linear estimate (\ref{linEst}) we know that $\nX{S(t)u_0}{X_1^{s,b}(B)}\leq c \nH{u_0}{s}$. From the definition of $X_T^{s,b}$ spaces we know that $T_1<T_2$ implies $X_{T_1}^{s,b}\subset X_{T_2}^{s,b}$. Therefore, for some $T<1$, $\nX{S(t)u_0}{X_T^{s,b}(B)}\leq c_0 \nH{u_0}{s}$.

Let us define $R'=2c_0R$.

From estimate (\ref{nlTerm}) we obtain, for $T<1$,
$$\nXsbT{\Phi(u)}{s}{b} \leq c_0\nH{u_0}{s} + c_1T^{1-b-b'}\nXsbT{u \bar{u} u}{s}{-b'}.$$

This, combined with (\ref{FNLest}), gives
$$\nXsbT{\Phi(u)}{s}{b} \leq c_0\nH{u_0}{s} + c_2T^{1-b-b'}\nXsbT{u}{s}{b}.$$
Taking $T<1$ such that $T^{1-b-b'} c_2 R'^3\leq c_0 R$, we ensure that $\Phi : B(0,R') \subset X_T^{s,b}(B) \rightarrow B(0,R') \subset X_T^{s,b}(B)$. In order to show that $\Phi$ is a contraction, let $u_1, u_2 \in B(0,R') \subset X_T^{s,b}(B)$. Then (\ref{linEst}) gives
$$\nXsbT{\Phi(u)-\Phi(u_2)}{s}{b} \leq c_2T^{1-b-b'}\nXsbT{|u_1|^2 u_1 - |u_2|^2 u_2}{s}{b}.$$
Decomposing $|u_1|^2 u_1-|u_2|^2 u_2 = u_1^2(\bar{u_1} - \bar{u_2}) + \bar{u_2}(u_1-u_2)(\bar{u_1}-\bar{u_2})$ and using (\ref{FNLest}), we obtain
$$\nXsbT{\Phi(u_1)-\Phi(u_2)}{s}{b} \leq c_3T^{1-b-b'}R'^2\nXsbT{u_1- u_2}{s}{b}.$$
Choosing $T<1$ eventually smaller and using (\ref{star}), we ensure that $\Phi$ is a contraction. The existence and uniqueness of $u\in X_T^{s,b}(B)$ such that $\Phi(u)=u$ follows. Since $b>\frac 12$, $u\in C((-T,T),H_D^s(B))$.

Let us show the Lipschitz property of the flow $u_0 \in B(0,R)\subset H_D^s(B) \rightarrow u \in X_T^{s,b}(B)$. Consider $u,v$ two solutions of (\ref{NLS}) with initial data $u_0,v_0 \in B(0,R)\subset H_D^s(B)$. From the Duhamel formula, $u$ and $v$ satisfy (\ref{2star}). As above, we show
$$\nXsbT{u-v}{s}{b} \leq c_0 \nH{u_0-v_0}{s} +c_3T^{1-b-b'}R'^2\nXsbT{u- v}{s}{b}.$$
Since $T$ was chosen such that $c_3T^{1-b-b'}R'^2<1$ we obtain $$\nXsbT{u-v}{s}{b}\leq c \nH{u_0-v_0}{s}.$$
\end{proof}

\subsection{Proof of the fundamental Lemma \ref{Flem}}
The fundamental Lemma \ref{Flem} helps us to control the cubic nonlinearity in a suitable $X^{s,b}(\RRB)$ space (imposed by estimate (\ref{nlTermEst})). We use a decomposition of the spectrum of functions $u_j \in X^{s,b}(\RRB)$. Obtaining suitable estimates helps us to sum over all frequencies by means of geometric sums.

A simple duality argument leads to the following equivalence : $u\in X^{s,b}(\RRB)$ $ \iff $ for all $u_0 \in X^{\infty, \infty}(\RRB) = \cap_{s>0, b\in \RR} X^{s,b}(\RRB)$ we have $$|<u,u_0>| \leq c \nXsb{u_0}{-s}{-b},$$ where $<,>$ denotes the duality bracket between $\mathcal{S'}$ and $\mathcal{S}$. Thus, by duality, (\ref{FNLest}) is implied by 
\begin{equation}
\label{dualF}
|\int_\RR \int_B u_0 u_1 u_2 u_3 {\rm d}x {\rm d}t| \leq c \prod_{j=1}^3 \nXsb{u_j}{s}{b} \nXsb{u_0}{-s}{b'},
\end{equation}
holding for all $u_0\in X^{\infty,\infty}(\RR \times B).$ We prove a similar result for spectrally localized functions and then sum over all frequencies. 

For $j\in \{0,1,2,3\}$, let $N_j \in \delta^{\mathbb{N}}$ (see Remark \ref{deltaNo}). We denote by $u_{jN_j}=\mathbbm{1}_{\sqrt{-\tr} \in {[N_j, \delta N_j]}} u_j$. Using the definition of $X^{s,b}(\RRB)$ spaces the following equivalence holds

\begin{equation}
\label{nXsbTrunc}
\nXsb{u_j}{s}{b}^2 \cong \sum_{N_j \in \delta^\mathbb{N}} \nXsb{u_{jN_j}}{s}{b}^2 \cong \sum_{N_j \in \delta^\mathbb{N}} N_j ^{2s} \nXsb{u_{jN_j}}{0}{b}^2.
\end{equation}

\noindent We denote by $\underline{N}=(N_0,N_1,N_2,N_3)$ the quadruple of $\delta^p$ numbers, $p\in \mathbb{N}$, and by $$I(\underline{N})=\int_{\RR\times B} \prod_{j=0}^3 u_{jN_j} {\rm d}t{\rm d}x.$$ 

In the following lemma we estimate $I(\underline{N})$ in two different ways. The first estimate is a consequence of property $\Ps$ and will be used in the case $N_0\leq c \max(N_1,N_2,N_3)$. The second one follows from an integration by parts and property $\Gs$. It will be useful in the case $N_0 \geq c \max(N_1,N_2,N_3)$. These two estimates imply Lemma \ref{Flem}.

\begin{lem}
\label{IntermedLem}
Under the assumptions $\Psz$ and $\Gsz$ for $s_0>0$, for all $s'>s_0$ there exists $0<b'<\frac 12$ , $c>0$ 
such that, assuming $N_3\leq N_2 \leq N_1$, the following inequalities hold :
\begin{equation}
\label{bLin}
|I(\underline{N})| \leq c (N_2 N_3)^{s'} \prod_{j=0}^3 \nXsb{u_{jN_j}}{0}{b'},
\end{equation}
\begin{equation}
\label{bLinIPP}
|I(\underline{N})| \leq c \left(\frac{N_1}{N_0}\right)^{2} (N_2 N_3)^{s'} \prod_{j=0}^3 \nXsb{u_{jN_j}}{0}{b'}
\end{equation}
\end{lem}

\begin{rem}
Please note that the assumption $N_3\leq N_2 \leq N_1$ does not restrict the generality since the estimate we prove is symmetric with respect to $u_1,u_2,u_3$.
\end{rem}

\begin{proof}
We start by proving inequality (\ref{bLin}). Using Holder inequality and (\ref{X014}), we obtain
\begin{eqnarray}
\label{R1}
|I(\underline{N})| &\leq & c \nLi{u_{3N_3}}{4}{L_x^\infty} \nLi{u_{2N_2}}{4}{L_x^\infty} \nLi{u_{1N_1}}{4}{L_x^2} \nLi{u_{0N_0}}{4}{L_x^2} \nonumber\\
&\leq & c(N_2 N_3)^ {\frac 32 + \epsilon} \prod_{j=0}^3 \nLi{u_{jN_j}}{4}{L_x^2} \leq c(N_2 N_3)^ {\frac 32} \prod_{j=0}^3 \nXsb{u_{jN_j}}{0}{\frac 14}.
\end{eqnarray}
We have used the Sobolev embedding $\nLi{u_N}{\infty}{B} \leq cN^{\frac 32 + \epsilon} \nLdi{u_N}{B}$. The inequality without $\epsilon$ also holds, but we don't need it here.

Using the Cauchy Schwarz inequality and (\ref{f2ISB}) (as we assumed $\Ps$), we obtain that for any $b_0>\frac 12$ there exist $c_0>0$ such that
\begin{eqnarray}
\label{R2}
|I(\underline{N})| &\leq &\nLdi{u_{0N_0} u_{2N_2}}{\RR\times B} \nLdi{u_{1N_1} u_{3N_3}}{\RR\times B} \nonumber \\
&\leq& c_0 (N_2 N_3)^{s_0} \prod_{j=0}^3 \nXsb{u_{jN_j}}{0}{b_0}.
\end{eqnarray}

In order to interpolate between (\ref{R1}) and (\ref{R2}), we use the decomposition $u_{jN_j}= \sum_{K_j} u_{jN_j K_j}$, where $u_{jN_j K_j} = \mathbbm{1}_{K_j \leq \langle i\partial_t + \tr \rangle <\delta K_j} u_{jN_j}$ and the sum is taken over $\delta^p$ numbers, for $p\in \mathbb{N}$ : $K_j \in \delta^\mathbb{N}$. 
Let us denote by $I(\underline{N},\underline{K})=\int_{\RR\times B} \prod_{j=0}^3 u_{jN_jK_j} {\rm d}x{\rm d}t.$
Estimates (\ref{R1}) and (\ref{R2}) give, under these settings, $$I(\underline{N},\underline{K}) \leq c (N_2 N_3)^\sigma (\prod_{j=0}^3 K_j)^\beta \prod_{j=0}^3 \nLdi{u_{jN_jK_j}}{\RRB},$$
where $(\sigma,\beta)$ equals either $(\frac 32+\epsilon,\frac 14)$ or $(s_0, b_0)$. 
For $s_0 <s' <1$ we can choose $\epsilon >0$, $b_0 >\frac 12$ and $0<b_1 < \frac 12$ such that by interpolation we have the same estimate for $(\sigma, \beta)= (s', b_1)$. 

Taking $b'\in (b_1,\frac 12)$, this reads
$$|I(\underline{N},\underline{K})| \leq c (N_2 N_3)^{s'} \prod_{j=0}^3 K_j^{b_1-b'} \nXsb{u_{jN_jK_j}}{0}{b'}.$$
Summing up over $\underline K \in (\delta^\mathbb{N})^4$, 
we obtain, by means of Cauchy Schwarz and geometric series ($\delta>1$) that
$$|I(\underline{N})| \leq c (N_2 N_3)^{s'} \prod_{j=0}^3 \nXsb{u_{jN_j}}{0}{b'},$$
which concludes the proof of (\ref{bLin}).

For the proof of (\ref{bLinIPP}) we use the Green formula : $$\int_B \tr f g - f \tr g {\rm d}x= \int_{\mathbb{S}^2} \frac{\partial f}{\partial \nu} g - f \frac{\partial g}{\partial \nu} {\rm d}\sigma.$$
Please note that $u_{0N_0} = \sum_{\mu_k \sim{} N_0} c_k e_k$, where $c_k=(u_{0N_0},e_k)_{L_2}$ and $e_k$ are eigenfunctions of the Dirichlet Laplacian associated with eigenvalues $\mu_k^2$. Using that $-\tr e_k = \mu_k^2 e_k$, we write 
$$u_{0N_0} = - \frac{\tr}{(N_0)^2} \sum_{\mu_k \sim{} N_0} c_k \left ( \frac{N_0}{\mu_k}\right)^2 e_k.$$ We define $Tu_{0N_0} =\sum_{\mu_k \sim{} N_0} c_k \left ( \frac{N_0}{\mu_k}\right)^2 e_k$ and $Vu_{0N_0} = \sum_{\mu_k \sim{} N_0} c_k \left ( \frac{\mu_k}{N_0}\right)^2 e_k$. Obviously $TVu_{0N_0} = VTu_{0N_0}= u_{0N_0}$ and $\nH{Tu_{0N_0}}{s} \sim \nH{u_{0N_0}}{s}$ for all $s$. Using this notation, $$u_{0N_0} = -\frac{\tr}{(N_0)^2} Tu_{0N_0}.$$ Placing this identity in $I(\underline{N})$ and using that ${u_{jN_j}} _{|_{\mathbb{S}^2}}=0$, we obtain by the Green formula that $$I(\underline{N})= \frac {1}{N_0^2} \int_{\RRB} Tu_{0N_0} \tr(u_{1N_1} u_{2N_2} u_{3N_3}).$$ The Laplacian distributes as follows : $\tr(u_{1N_1} u_{2N_2} u_{3N_3})$ equals
$$\sum (\tr u_{\sigma_1N_{\sigma_1}}) u_{\sigma_2 N_{\sigma_2}} u_{\sigma_3 N_{\sigma_3}} + \sum (\triangledown u_{\sigma_1 N_{\sigma_1}})\cdot (\triangledown u_{\sigma_2 N_{\sigma_2}}) u_{\sigma_3 N_{\sigma_3}},$$ where the sum is taken over all permutations $\sigma \in \Sigma_3$. We consider a representative term of each sum. Let $$J_{11}(\underline{N})=  \int_{\RRB} Tu_{0N_0}  \tr u_{1N_1} u_{2N_2} u_{3N_3}$$ and $$J_{12}(\underline N) = \int_{\RRB} Tu_{0N_0} (\triangledown u_{1N_1})\cdot (\triangledown u_{2N_2}) u_{3N_3}.$$ We then have that $I(\underline{N})$ is a sum of terms similar to $\frac{1}{N_0^2}J_{11}(\underline{N})$ or to $\frac{1}{N_0^2}J_{12}(\underline{N})$. As we will see later on, those are the largest terms in each sum. Using that $\tr u_{N} = -N^2 Vu_N$, $J_{11}(\underline{N}) = -N_1^2 \int_{\RRB} Tu_{0N_0}  Vu_{1N_1} u_{2N_2} u_{3N_3}$.
We can estimate it by (\ref{bLin}) and use the norm equivalence $\nH{Tu_{0,N_0}}{s} \sim \nH{u_{0,N_0}}{s} \sim \nH{Vu_{0,N_0}}{s}$ to conclude that the first kind of terms behave
 as announced : $$\frac{1}{N_0^2}|J_{11}(\underline{N})| \leq c\frac{N_1^2}{N_0^2} (N_2 N_3)^{s'} \prod_{j=0}^3 \nXsb{u_{jN_j}}{0}{b'}.$$ 

The bounds on $J_{12}(\underline{N})$ will follow from a similar approach to the one used to prove (\ref{bLin}). Let us recall what we did in order to bound $I(\underline{N})$ in (\ref{bLin}). In order to control $I(\underline{N})$ we have obtained two estimate, (\ref{R1}) and (\ref{R2}), and we have interpolated between them using a supplementary decomposition. 

Estimate (\ref{R1}) was obtained via Holder inequality, estimating the two functions localized at small frequencies in $L_t^4L_x^\infty$ and the other two in $L_t^4L_x^2$. A loss of $(N_2N_3)^{\frac 32+\epsilon}$ arose from changing the $L^\infty(B)$ norm into a $L^2(B)$ norm. Here a supplementary loss of $N_1 N_2$ comes from estimating the gradient of a localized function in $L^2$ or $L^\infty$ norm : $\nLdi{\nabla u_{N}}{B} \leq c N \nLdi{u_{N}}{B}$ and by the Sobolev embedding $\nLi{\nabla u_N}{\infty}{B} \leq c N^{\frac 52 +\epsilon}.$ Therefore, the estimate corresponding to (\ref{R1}) reads :
$$|J_{12}(\underline N)| \leq c (N_1 N_2) (N_2N_3)^{\frac 32 + \epsilon} \prod_{j=0}^3 \nXsb{u_{jN_j}}{0}{\frac 14}.$$

Estimate (\ref{R2}) was obtained using the Cauchy Schwarz inequality, coupling a function localized at high frequency with a function localized at a smaller frequency. Here we have to do the same and moreover to do so without putting two gradients together :
$$|J_{12}(\underline N)| \leq c \nLdi{Tu_{0N_0} \nabla u_{2N_2}}{\RRB} \nLdi{\nabla u_{1N_1} u_{3N_3}}{\RRB}.$$
Using (\ref{f3ISB}) (as we assume $\Gsz$), we obtain, just as in (\ref{R2}), that :
$$|J_{12}(\underline N)| \leq c (N_1 N_2) (N_2N_3)^{s_0} \prod_{j=0}^3 \nXsb{u_{jN_j}}{0}{b_0}.$$

Please note that compared with (\ref{R1}) and (\ref{R2}), the estimates on $J_{12}(\underline N)$ have just a supplementary factor of $N_1 N_2$. Therefore the interpolation argument works without a change and we obtain
$$\frac{1}{N_0^2}|J_{12}(\underline{N})| \leq c\frac{N_1 N_2}{N_0^2} (N_2 N_3)^{s'} \prod_{j=0}^3 \nXsb{u_{jN_j}}{0}{b'}.$$
Please note that $N_1 N_2$ is the worst term that comes out in estimating $J_{12}(\underline{N})$ (under the assumption $N_3 \leq N_2 \leq N_1$) and it is smaller than the $N_1^2$ that arose in the estimate of $J_{11}(\underline{N})$.
\end{proof}

Let us show how Lemma \ref{IntermedLem} implies Lemma \ref{Flem}.
\begin{proof} (of Lemma \ref{Flem}) We want to bound $I=\int_{\RRB} u_0 u_1 u_2 u_3 $ as announced in (\ref{dualF}). We shall decompose $I= \sum_{\underline{N}\in  (\delta^{\mathbb{N}})^4} I(\underline{N})$, where $$I(\underline{N})=\int_{\RRB} u_{0N_0} u_{1N_1} u_{2N_2} u_{3N_3}.$$

By symmetry we can consider that $N_3\leq N_2 \leq N_1$. We separate the sum into two cases : $\sum\prime=\sum_{\underline{N} : N_0 \leq c N_1}$ and $\sum{\prime\prime}=\sum_{\underline{N} : N_0>c N_1}$.

Let $s'$ be such that $\frac 12 <s' <s$. Using Lemma \ref{IntermedLem} and (\ref{nXsbTrunc}), we have 
$$|\sum\prime I(\underline{N})| \leq c \sum\prime(N_2 N_3)^{s'-s} \left( \frac{N_0}{N_1}\right)^s \nXsb{u_{0N_0}}{s}{-b}\prod_{j=1}^3 \nXsb{u_{jN_j}}{s}{b'}.$$ Using the Cauchy Schwarz inequality and the equivalence(\ref{nXsbTrunc}) we bound $$|\sum\prime I(\underline{N})| \leq c \nXsb{u_2}{s}{b'} \nXsb{u_3}{s}{b'} \sum_{N_0 \leq cN_1} \left( \frac{N_0}{N_1}\right)^s \alpha(N_0) \beta(N_1),$$ where we denote by $\alpha(N_0)=\nXsb{u_{0N_0}}{-s}{b'}$ and by $\beta(N_1)=\nXsb{u_{1N_1}}{s}{b'}$. From (\ref{nXsbTrunc}) we have $$\sum_{N_0} \alpha(N_0)^2 \cong \nXsb{u_0}{-s}{b'}^2\ {\rm and}\ \sum_{N_1} \beta(N_1)^2 \cong \nXsb{u_1}{s}{b'}^2.$$ Since both $N_0$ and $N_1$ are $\delta^p$ numbers, $p\in \mathbb{N}$, we can write $N_1 = \delta^l N_0$, for $N_0\geq N(l)=\max(1,\delta^{-l})$, where $l$ is an integer, $l\geq -l_0$, for some $l_0 \in \mathbb{N}$, depending only on $c$. Thus,
\begin{eqnarray*}
\sum_{N_0 \leq cN_1} \left( \frac{N_0}{N_1}\right)^s \alpha(N_0) \beta(N_1) = \sum_{l\geq -l_0} \sum_{N_0\geq N(l)} \delta^{-sl} \alpha(N_0) \beta(\delta^l N_0)\\
\leq  c\sum_{l>-l_0} \delta^{-sl} \left( \sum_{N_0} \alpha(N_0)^2 \right)^{\frac 12} \left( \sum_{N_0 \geq N(l)} \beta(\delta^l N_0)^2 \right)^{\frac 12}\\ 
\leq c\nXsb{u_0}{-s}{b'} \nXsb{u_1}{s}{b'}.
\end{eqnarray*}
Noticing that $\nXsb{u}{s}{b'} \leq \nXsb{u}{s}{b}$ for $b>b'$, we have
$|\sum\prime I(\underline{N})| \leq \prod_{j=1}^3 \nXsb{u_j}{s}{b} \nXsb{u_0}{-s}{b'}.$

We bound the second sum in a similar manner, using estimate (\ref{bLinIPP}) of Lemma \ref{IntermedLem} : $$|\sum{\prime\prime} I(\underline{N})| \leq c \sum{\prime\prime}(N_2 N_3)^{s'-s}  \left( \frac{N_1}{N_0}\right)^{2-s} \nXsb{u_{0N_0}}{-s}{b'}\prod_{j=1}^3 \nXsb{u_{jN_j}}{s}{b'}.$$
Since $N_0\geq cN_1$ we write $N_0=\delta^p N_1$, for $p\in \mathbb{N}$, $p>p_0$ depending only on $c$. 
As above, inverting the role of $N_0$ and $N_1$, we obtain $$|\sum{\prime\prime} I(\underline{N})| \leq c \nXsb{u_0}{-s}{b'} \nXsb{u_1}{s}{b'} \nXsb{u_2}{s}{b'} \nXsb{u_3}{s}{b'}.$$
\end{proof}

\section{Bilinear Strichartz estimates : proof of Proposition \ref{bilinLem}}
\label{bilinear}
We give the proof of those estimate in the case of radial data on $B$, the unit ball of $\RR^3$. The proof makes use of the exact form of the eigenfunctions of the Dirichlet Laplacian in $L_{rad}^2(B)$.

\subsection{Eigenfunctions of the Dirichlet Laplacian}
\label{EDL}
This subject is covered by almost any book of mathematical methods for physicists. 
 We remind them here mostly for notation reasons. 

The eigenfunctions of the Dirichlet Laplacian in $L^2(B)$, where $B\subset \RR^3$, are of the form
$$\varphi_{n,k}(x) = c|x|^{-\frac 12} J_{n+\frac 12} (\jnk |x|) H_n\left( \frac {x}{|x|} \right),$$ where $J_\nu$ denotes a Bessel function of order $\nu$ (for more details on $J_\nu$ see Section \ref{BesselFunct}), $\jnk$ denotes the $k$th zero of the $J_{n+\frac 12}$ and $H_n$ is a spherical harmonics of order $n$. This $\varphi_{n,k}$ is associated to eigenvalues $\jnk ^2$ : $-\tr \varphi_{n,k} = \jnk^2 \varphi_{n,k}$. For more properties of Bessel functions see Section \ref{BesselFunct}.

The eigenfunctions of the Dirichlet Laplacian in $L_{rad}^2$ are the $\varphi_{n,k}$ corresponding to $n=0$. Let us denote by $r=|x|$ and by $ e_k(r)=\varphi_{0,k}(x)= c r^{-\frac 12} J_{\frac 12} (j_{\frac 12,k} r)$. Since $J_{\frac 12}(r) = c \frac {\sin (r)}{r^\frac 12}$, we get $e_k(r) =c \frac{\sin(k\pi r)}{r}$. The eigenfunction $e_k$ is associated to eigenvalue $k^2 \pi^2$. The constant $c$ is chosen such that $\nLdi{e_k}{B}=1$. Those eigenfunctions form an eigenbase of $L_{rad}^2(B)$.

The standard approach for proving properties $\Ps$ and $\Gs$ for $S(t)$ is to decompose $u_0$ and $v_0$ on the base $(e_k)$.

\subsection{Proof of properties $\Ps$ and $\Gs$ for radial data and $s>\frac 12$}

Let us remind how property $\Ps$ reads : for all $N,L>0$ and $u_0, v_0 \in L^2(B)$ satisfying $\mathbbm{1}_{\sqrt{-\tr} \in {[N,\delta N]}} u_0=u_0,\ \ \mathbbm{1}_{\sqrt{-\tr} \in {[L,\delta L]}} v_0=v_0,$
$$\nLdi{S(t)u_0 S(t)v_0}{(0,1)\times B} \leq C (\min(N,L))^{s} \nLdi{u_0}{B}\nLdi{v_0}{B}. $$
If moreover $u_0, v_0 \in H_0^1(B)$, $S(t)$ satisfies property $\Gs$ if
$$\nLdi{\triangledown S(t)u_0 S(t)v_0}{(0,1)\times B} \leq C N(\min(N,L))^{s} \nLdi{u_0}{B}\nLdi{v_0}{B}.$$

We write $u_0 (x) = \sum_{k \sim_\delta  N } c_k e_k(r)$ and $v_0 (x) =\sum_{l \sim_\delta  L } d_l e_l(r)$, where $c_k = (u_0,e_k)_{L^2}$ and $d_l = (v_0,e_l)_{L^2}$. Then $S(t)u_0 = \sum_{k \sim_\delta  N } e^{-itk^2 \pi^2 } c_k e_k(r)$ and $S(t)v_0 = \sum_{l \sim_\delta  L } e^{-itl^2 \pi^2 } d_l e_l(r)$. Thus, in order to prove $\Ps$, what we want to estimate is the $L^2((0,T)\times B)$ norm of $$E(N,L)=\sum_{k\sim_\delta  N,\ l \sim_\delta  L} e^{-it(k^2+l^2) \pi^2} c_k d_l e_k e_l.$$ For $\Gs$ we need to estimate $$E_1(N,L)=\sum_{k\sim_\delta  N,\ l \sim_\delta  L} e^{-it(k^2+l^2) \pi^2} c_k d_l (\nabla e_k) e_l.$$

We proceed to the proof of $\Ps$.
\subsubsection{Property $\Ps$}
\label{proofPsrad}

For $\tau \in \mathbb{N}$, let us define $$\Lambda_{NL}(\tau) = \{ (k,l)\in \mathbb{N}^2\ s.t. \ k\sim_\delta  N,\ l\sim_\delta  L,\ k^2 + l^2 = \tau \}.$$ Applying Parseval in time variable, we have $$\nLdi{E(N,L)}{(0,\frac{2}{\pi})\times B}^2 = c \sum_{\tau \in \mathbb{N}} \nLdi{\sum_{(k,l)\in \Lambda_{NL}(\tau)} c_k d_l e_k e_l}{B}^2.$$ By means of Cauchy Schwarz we obtain
\begin{equation}
\label{cBest}
\nLdi{E(N,L)}{(0,\frac{2}{\pi})\times B}^2 \leq c \sum_{\tau \in \mathbb{N}} \# 
\Lambda_{NL} (\tau) \sum_{(k,l)\in \Lambda_{NL}(\tau)} |c_k d_l|^2  \nLdi{e_k e_l}{B}^2,
\end{equation} 
where $\#\Lambda_{NL}(\tau)$ denotes the cardinal of the $\Lambda_{NL}(\tau)$ set. We present two lemmas that give us the bound on the elements of the right hand side in (\ref{cBest}). First we remind an estimate on the number of lattice points on an arc of a circle in $\RR^2$. 

\begin{citelem} (3.2 of \cite{BGTinv}) For all $\epsilon>0$, $\#\Lambda_{NL}(\tau) = O(\min(N,L)^\epsilon)$.
\end{citelem}
The following lemma gives the bilinear estimate of eigenfunctions in $L^2(B)$.
\begin{lem}
\label{Beigen}
 There exists $c>0$ such that for $k,l\in \mathbb{N}$, $$\nLdi{e_k e_l}{B} \leq c \min(k,l)^{\frac 12}.$$
\end{lem}
\begin{proof} 
Replacing $e_k$ by their exact expression, we obtain
$$\nLdi{e_k e_l}{B}^2 = c\int_0^1 \frac{\sin^2(k\pi r) \sin^2(l\pi r)}{r^2} {\rm d}r.$$ 
Let us consider $k\leq l$. This does not restrict the generality. We perform the change of variable $r=\frac {s}{k\pi}$ :
$$\nLdi{e_k e_l}{B}^2 = c k \int_0^{k\pi} \frac{\sin^2(s) \sin^2(\frac{l}{k}s)}{s^2} {\rm d}s \leq ck \int_0^\infty \frac{\sin^2(s)}{s^2} {\rm d}s .$$
As $\int_0^\infty \frac{\sin^2(s)}{s^2} {\rm d}s <\infty$, the result follows.
\end{proof}

Using the two previous lemmas and (\ref{cBest}), we have
$$\nLdi{E(N,L)}{(0,\frac{2}{\pi})\times B}^2 \leq c_\epsilon \min(N,L)^{1+\epsilon} \sum_{k\sim N,\ l\sim L} |c_k|^2 |d_l|^2.$$
Since $\sum_{k\sim N,\ l\sim L} |c_k|^2 |d_l|^2 = \nLdi{u_0}{B}^2 \nLdi{v_0}{B}^2$, we obtain property $\Ps$ on radial data, for $s>\frac 12 +\epsilon$.

\subsubsection{Property $\Gs$}
Just as we did for estimating $E(N,L)$ in (\ref{cBest}), applying Parseval in time variable and Cauchy Schwarz in $x$ variable, we get 
\begin{equation}
\label{cDBest}
\nLdi{E_1(N,L)}{(0,\frac{2}{\pi})\times B}^2 \leq c \sum_{\tau \in \mathbb{N}} \# 
\Lambda_{NL} (\tau) \sum_{(k,l)\in \Lambda_{NL}(\tau)} |c_k d_l|^2  \nLdi{(\nabla e_k) e_l}{B}^2,
\end{equation} 
The following lemma gives the crucial estimate for proving $\Gs$.
\begin{lem} 
\label{DBeigen}
There exists $c>0$ such that for $k,l\in \mathbb{N}$, 
\begin{equation}
\label{gradBeigen}
\nLdi{(\nabla e_k) e_l}{B} \leq c k \min(k,l)^{\frac 12}.
\end{equation}
\end{lem}
\begin{proof} 
Please note that for radial data $(\nabla_x f)\cdot (\nabla_x g) = (\partial_r f)(\partial_r g)$. Thus, $$\nLdi{(\nabla e_k) e_l}{B}^2 = c\int_0^1 \left ( \partial_r \frac {\sin(k\pi r)}{r} \right ) ^2 \sin^2(l\pi r){\rm d}r.$$
Expanding the derivative  $\left ( \partial_r \frac {\sin(k\pi r)}{r} \right )$ we obtain
\begin{equation}
\label{explicitGBE}
\left ( \partial_r \frac {\sin(k\pi r)}{r} \right ) ^2 = \frac{k^2 \pi^2}{r^2} \left (\cos (k\pi r)-  \frac {\sin(k\pi r)}{k \pi r} \right ) ^2.
\end{equation}

As opposed to Lemma \ref{Beigen}, $k$ and $l$ are not in a symmetric position. Therefore, we shall treat two cases : $ k\leq l$ and $l<k$.

In the first case, $ k\leq l$, (\ref{gradBeigen}) reads $\nLdi{(\nabla e_k) e_l}{B}^2 \leq c k^3.$ We perform the change of variable $r=\frac {s}{k\pi}$ in (\ref{explicitGBE}). Then 
$$\nLdi{(\nabla e_k) e_l}{B}^2 = k^3 \pi^3 \int_0^{k\pi}\frac{\sin^2\left(\frac{l}{k}s \right)}{s^2}\left(\cos(s)-\frac{\sin(s)}{s} \right)^2  ds.$$
We bound $\sin^2(\frac lk s)\leq 1$ and take the integral from $0$ to $\infty$. 
Close to $0$ we have $\left(\cos(s)-\frac{\sin(s)}{s} \right)^2 = \frac{s^4}{9} + O(s^6)$ and therefore the integral converges. Away from $0$ we bound $\left(\cos(s)-\frac{\sin(s)}{s} \right)^2 \leq c$ and $\frac{1}{s^2}$ provides convergence. Consequently $\int_0^\infty \frac {1}{s^2}\left(\cos(s)-\frac{\sin(s)}{s} \right)^2 ds \leq c$ and $\nLdi{(\nabla e_k) e_l}{B}^2 \leq c k^3.$

In the second case, $l<k$, (\ref{gradBeigen}) reads $\nLdi{(\nabla e_k) e_l}{B}^2 \leq c k^2 l.$ We perform the change of variable $r=\frac {s}{l\pi}$ in (\ref{explicitGBE}). Then 
$$\nLdi{(\nabla e_k) e_l}{B}^2 = k^2 l\pi^3 \int_0^{l\pi}\frac{\sin^2(s)}{s^2} \left(\cos\left( \frac{k}{l}s \right) - \frac{\sin\left(\frac{k}{l}s \right)}{\frac{k}{l}s} \right)^2  ds.$$ We have $\left(\cos(s) - \frac{\sin(s)}{s} \right)^2 \leq c$ for all $s\in \RR$. Consequently, $$\nLdi{(\nabla e_k) e_l}{B}^2 \leq c k^2 l \int_0^\infty \frac{\sin^2(s)}{s^2} ds \leq c k^2 l.$$ 
\end{proof}

\subsection{Property $\Ps$ in a particular non radial case}

The general method of obtaining local existence results from a bilinear estimate invites us to apply it to a general non radial data. One of the problems we encountered is the distribution of zeros of all Bessel functions of half integer order. 

Although we were not able to prove bilinear estimate for general data, we obtained encouraging results for data which are superpositions of eigenfunctions corresponding to the first zero of the Bessel function $J_{n+\frac 12}$. It is expected that this case should be the worst behaved. A similar analysis could be performed for $\nu=n$ (two space domains for example).

\subsubsection{Some properties of Bessel functions}
\label{BesselFunct}

Bessel functions are some of the most studied special function. They appear in many physical problems. The Bessel function $J_\nu$, of first kind of order $\nu$, satisfies the Bessel equation of order $\nu$ :
$$z^2 w{''} + z w{'} +(z^2-\nu^2)w=0.$$
If $\nu=n$ is an integer then $J_n$ are known as Bessel coefficients. If $\nu=n+\frac 12$ is half of an odd integer, then $J_\nu$ reduces to a finite combination of elementary functions (see e.g. \cite{Wat} 4.7 to 4.75).  Although we are interested only in Bessel functions of order $\nu=n+\frac 12$, we shall mainly use Schläfli's integral representation (see \cite{Wat}, p.176) : for $x\in \mathbb{R}_+$
\begin{equation}
\label{defJnu}
J_\nu (x)= \frac{1}{2\pi} \int_{-\pi}^\pi e^{i(x \sin t - \nu t)} {\rm d}t -\frac{\sin (\nu \pi)}{\pi} \int_0^\infty e^{-(x \sinh t +\nu t)}{\rm d}t.
\end{equation}
Let us denote the first term of the sum in (\ref{defJnu}) by $$T_{1,\nu}(x)= \frac{1}{2\pi} \int_{-\pi}^\pi e^{i(x \sin t - \nu t)} {\rm d}t$$  and the second term by $$T_{2,\nu}(x)=\frac{\sin (\nu \pi)}{\pi} \int_0^\infty e^{-(x \sinh t +\nu t)}{\rm d}t.$$ Thus, $J_\nu (x)=T_{1,\nu}(x)-T_{2,\nu}(x)$. 

Please note that for $\nu=n\in \mathbb{N}$, $T_{2,\nu}(x)$ vanishes and $J_n(x)=T_{1,n}(x)$. This expression was used in \cite{BGTgafa} to obtain asymptotics of the Bessel function $J_n$. However, their analysis does not use the fact that $\nu=n$ is an integer and can be performed for $\nu>0$. In the sequel we shall take  $\nu\in \frac 12 \mathbb{N}$.

The asymptotic behavior of Bessel function $J_\nu(x)$ is different according to the order $\nu$ 	and the variable $x$. Since they can get both very large, we distinguish the case of large order, of large variable and the transitional regions. As the coupling $\frac x\nu$ counts, we write $x=\nu \rho$, $\rho\in \RR_+$. We shall refer to \cite{BGTgafa} for the proof of the following two asymptotics of $T_{1,\nu}(x)$. 

\begin {citelem} Let $\beta_1 >0$ and $\gamma \in (0,\frac 12)$.
\begin{itemize}
	\item {Then for $\rho \in [1-\beta_1\nu^{-\gamma} , 1+\beta_1 \nu^{-\gamma}]$, $\nu \gg 1$ one has
	\begin{equation}
	\label{AiryT1nu}
	T_{1,\nu}(\nu \rho) = \left( \frac{2}{\nu \rho}\right)^{\frac 13} Ai(\nu^{\frac 23} (1-\rho)(\frac 2 \rho) ^\frac 13) + r_\nu(\rho),
	\end{equation}
	where $|r_\nu(\rho)| \leq c\nu ^{-\gamma_1}$, $\gamma_1<3\gamma -1$.}
	
	\item{ Then for any $k>0$ there exists $c_k>0$ such that for $\rho \in [0,1-\beta\nu^{-\gamma}]$, $\nu \gg 1$ one has :
	\begin{equation}
	\label{decreaseT1nu}
	|T_{1,\nu}(\nu \rho) | \leq c_k \nu ^{-k}.
	\end{equation}}
\end{itemize}
\end{citelem}

The second term of the sum in (\ref{defJnu}), $T_{2,\nu}(\nu \rho)=\frac{\sin (\nu \pi)}{\pi} \int_0^\infty e^{-\nu(\rho \sinh t +t)}{\rm d}t$,  can be developed in terms of $\nu$ and $\rho$, but we shall only use the simple remark that there exists $c>0$ such that for all $\rho>0$,
\begin{equation}
\label{decreaseT2nu}
|T_{2,\nu}(\nu \rho)| \leq c \nu^{-1}.
\end{equation}

The Bessel functions have infinitely many zeros. In section \ref{EDL} we have denoted by $j_{\nu,k}$ the $k$th zero of the Bessel function $J_\nu$. We denote by $z_\nu = j_{\nu,1}$, the first positif zero of $J_{\nu}$. It is known that in the neighborhood of $z_\nu$ the Bessel function $J_{\nu}$ behaves like an Airy function (\ref{AiryT1nu}). The asymptotics of the first zero of the Bessel function reads as follows

\begin{lem}
\label{asymptzn1}
There exists $\alpha>0$ and $\beta$ two real constants such that, for all $\nu \in \RR_+$, $$z_{\nu} = \nu + \alpha \nu^{\frac 13} + \beta \nu^{-\frac 13} + O(\nu^{-\frac 23}).$$ 
\end{lem}
The asymptotics of the first zero of the Bessel function has been studied for example by the Sturm method (see \cite{Wat}, ch. XIV). This gives $z_{\nu} = \nu + \alpha \nu^{\frac 13} + O(\nu^{-\frac 13}).$ For reasons we shall present further, we need more precise asymptotics. For completeness and for the convenience of the reader, we shall give the proof in Appendix \ref{pfAsymptzn1}.

Combining these asymptotics, we obtain the following bilinear estimate
\begin{lem}
\label{bilinJn1Jp1}
 For $\nu ,p \in \frac 12 \mathbb{N}$ such that $\nu \sim N$ and $p\sim L$, the following estimate holds :
\begin{equation}
\label{estJnJp}
\nLdi{r^{-1} J_{\nu}(z_\nu r) J_{p}(z_{p} r)}{B} \leq c \min(N,L)^\frac 13 \nLdi{r^{-\frac 12} J_{\nu} (z_\nu r)}{B} \nLdi{r^{-\frac 12} J_{p} (z_p r)}{B}.
\end{equation}
\end{lem}

\begin{proof}We mention that for this proof, the precise asymptotics from Lemma \ref{asymptzn1} does not provide supplementary informations with respect to the use of Sturm asymptotics. From Lemma 5.1 of \cite{BGTgafa} we know that $\nLdi{r^{-\frac 12} J_{\nu} (z_\nu r)}{B}^2 \approx \nu^{-\frac 43}$. In order to estimate the left hand side term in (\ref{estJnJp}), we follow the same approach as in the proof of Lemma 5.1 from \cite{BGTgafa}.
We want to estimate 
\begin{equation}
\label{idJnJp}
\nLdi{r^{-1} J_{\nu}(z_\nu r) J_{p}(z_{p} r)}{B}^2 = c \int_0^1 |J_{\nu}(z_\nu r) J_{p}(z_{p} r)|^2 {\rm d}r.
\end{equation}
Let us consider $N<L$ and $\nu <p$. With the notations of the lemma on the asymptotics of $T_{1,\nu}$, let us denote by $r_{\nu}=z_\nu (1-\beta_1  \nu^{-\gamma})$. As $\beta_1>0$ and $\nu<p$ then $$r_{\nu}< r_p.$$ Following the asymptotics of $J_\nu$ and $J_p$, we decompose the integral $\int_0^1 = \int_0^{r_\nu} + \int_{r_\nu}^{r_p} + \int_{r_p}^1$. Combining (\ref{decreaseT1nu}) and (\ref{decreaseT2nu}), we know that for $0<r<r_\nu<r_p$, $J_\nu(z_\nu r)$ and $J_p(z_p r)$ decrease respectively like $\nu^{-1}$ and $p^{-1}$.  We deduce that 
$$\int_0^{r_\nu} |J_{\nu}(z_\nu r) J_{p}(z_{p} r)|^2 {\rm d}r \leq c \nu^{-2} p^{-2}.$$

We present the method for estimating $\int_{r_p}^1 |J_{\nu}(z_\nu r) J_{p}(z_{p} r)|^2 {\rm d}r.$  
For $r_\nu<r_p<r<1$ we know from (\ref{AiryT1nu}) that the asymptotic expansions of $J_\nu(z_\nu r)$ and $J_p(z_p r)$ are close to the Airy function. We bound the error terms by the main terms and perform the change of variable $r=\frac{p}{z_p}\rho$. Consequently,
$$\int_{r_p}^1 |J_{\nu}(z_\nu r) J_{p}(z_{p} r)|^2 {\rm d}r \leq \left( \frac{c}{z_\nu p} \right ) ^{\frac 23}\int_{1-\beta p^{-\gamma}}^{1+\alpha p^{-\frac 23}} Ai^2(\nu^{\frac 23} a(\frac{z_\nu}{\nu} \frac {p}{z_p}\rho)) Ai^2(p^{\frac 23} a(\rho)) \rho^{-\frac 43} {\rm d}\rho,$$
where we have denoted by $a(\rho)=(1-\rho)\left( \frac 2 \rho \right)^{\frac 13}$. The change of variable $\rho = 1+p^{-\frac 23}y$ yields to
$$\int_{r_p}^1 |J_{\nu}(z_\nu r) J_{p}(z_{p} r)|^2 {\rm d}r \leq \left( \frac{c}{\nu p^2} \right ) ^{\frac 23} c(n,p),$$
where $c(n,p)\rightarrow \int_{-\infty}^\alpha Ai^2(2^\frac 13 y) Ai^2(\alpha+ c_0(y-\alpha))dy$ if $\frac np \rightarrow c_0<1$ and $c(n,p) \rightarrow \int_{-\infty}^\alpha Ai^2(2^\frac 13 y) Ai^2(\alpha)dy$ if $\frac np \rightarrow 0$.

The estimate $\int_{r_\nu}^{r_p}|J_{\nu}(z_\nu r) J_{p}(z_{p} r)|^2 {\rm d}r \leq c \nu^{-\frac 43} p^{-2}$ is obtained via a similar computation, using that, for $r_\nu <r<r_p$, $J_p(z_p r)$ decreases like $p^{-1}$ (see (\ref{decreaseT1nu}) and (\ref{decreaseT2nu})) and $J_\nu(z_\nu r)$ behaves like the Airy function (\ref{AiryT1nu}). Coupling the estimate of each term of the decomposition $\int_0^1 = \int_0^{r_\nu} + \int_{r_\nu}^{r_p} + \int_{r_p}^1$ gives $\int_0^1 |J_{\nu}(z_\nu r) J_{p}(z_{p} r)|^2 {\rm d}r \leq c \nu ^{-\frac 23} p^{-\frac 43}$. Combined with $\nLdi{r^{-\frac 12} J_{\nu} (z_\nu r)}{B}^2 \approx \nu^{-\frac 43}$, this yields the result.
\end{proof}

\subsubsection{Proof of property $\Ps$ for data superposition of eigenfunctions corresponding to the first zero of $J_{n+\frac 12}$}
\label{PsZn1}
Let us denote by $\mathcal{V}$ the vector space spanned by the eigenfunctions corresponding to the first zero of the Bessel function. With the notations of section \ref{EDL}, $\mathcal{V} = \{u\in L^2(B) \ s.t. \ u=\sum_{n\in T} c_n \varphi_{n,1},\ c_n\in\mathbb{C}, \ \#T<\infty \}.$
Let us consider initial data $u_{0N}, v_{0L} \in\mathcal{V}$ localised at frequencies $N,L \in \delta^\mathbb{N}$ (see Remark \ref{deltaNo}) :

$$u_{0N} = \sum_{{z_{n+\frac 12}} \sim N} c_n \varphi_{n,1}\ and \ u_{0L} = \sum_{{z_{p+\frac 12}} \sim L} d_p \varphi_{p,1}.$$ Then 
\begin{equation}
\label{bilinSJnSJp}
S(t)u_{0N} S(t)v_{0L}=\sum_{{z_{n+\frac 12}} \sim N, {z_{p+\frac 12}} \sim L} e^{-it (z_{n+\frac 12}^2 + z_{p+\frac 12}^2)} c_n d_p \varphi_{n,1}(x)\varphi_{p,1}(x).
\end{equation}
We follow the same approach as in the radial case (see section \ref{proofPsrad}). The following general lemma allows us to handle the integral in time.

\begin{citelem}(5.2 of \cite{BGTens}) Let $\Lambda$ be a countable set of real numbers. Then for every $T>0$ there exists $c_T>0$ such that for every sequence $(a_\lambda)$ indexed by $\Lambda$, one has
$$\nLdi{\sum_{\lambda \in \Lambda} a_\lambda e^{i\lambda t}}{0,T} \leq c_T \left ( \sum_{l\in \mathbb{Z}} (\sum_{\lambda : |\lambda -l|\leq \frac 12} |a_\lambda|)^2 \right)^{\frac 12}.$$
\end{citelem}

Denoting by 
\begin{equation}
\label{defGNL}
\Gamma_{NL}(\tau) = \left\{ (n,p)\in \mathbb{N}^2\ s.t. \ z_{n+\frac 12}\sim N,\ z_{p+\frac 12}\sim L,\ |z_{n+\frac 12}^2 + z_{p+\frac 12}^2 - \tau| \leq \frac 12 \right\}
\end{equation} and using Cauchy Schwarz inequality we obtain $\nLdi{S(t)u_{0N} S(t)v_{0L}}{(0,1)\times B}^2$ bounded by
\begin{equation}
\label{mainPsJn1}
\sum_{\tau \in \mathbb{N}} \#\Gamma_{NL}(\tau) \sum_{(n,p) \in \Gamma_{NL}(\tau)} |c_n d_p|^2 
\int_0^1 |J_{n+\frac 12} (z_{n+\frac 12} r) J_{p+\frac 12}(z_{p+\frac 12} r)|^2 dr \nLdi{H_n H_p}{\mathbb{S}^2}^2.
\end{equation}

From Lemma \ref{bilinJn1Jp1} we know the behavior of the integral in $r$. In order to estimate the spherical harmonics we use a result by N.Burq, P.Gérard and N.Tzvetkov \cite{BGTinv} (see also \cite{BGTens}).
\begin{citelem} There exists $C>0$ such that, if $H_n$ and $H_p$ are two spherical harmonics of respective degrees $n$ and $p$,
\begin{equation}
\label{spharmonics}
\nLdi{H_n H_p}{\mathbb{S}^2} \leq c \min(n,p)^\frac 14 \nLdi{H_n}{\mathbb{S}^2} \nLdi{H_p}{\mathbb{S}^2}.
\end{equation}
\end{citelem}

We have to estimate the cardinal of $\Gamma_{NL}(\tau)$ in order to obtain a bilinear estimate of $S(t)u_{0N} S(t)v_{0L}$ from (\ref{bilinSJnSJp}). The following lemma is essential in the proof of Proposition \ref{ISBzn1}.
\begin{lem}
\label{cardGnl}
For $N,L \in \delta^\mathbb{N}$ we have $\#\Gamma_{NL}(\tau) \leq c \min(N,L)^{\frac 23}$, where $\Gamma_{NL}(\tau)$ is defined by (\ref{defGNL}).
\end{lem}
The proof of this lemma uses classical techniques in analytic number theory. We present it in Appendix \ref{proofCardGnl}. Let us show how we  deduce property $\Ps$ from the previous lemmas.

\begin{proof} (of Proposition \ref{ISBzn1}) Replacing into (\ref{mainPsJn1}) the results of Lemmas \ref{bilinJn1Jp1} and \ref{cardGnl}, as well as the bilinear estimate of the spherical harmonics in $L^2(B)$ (\ref{spharmonics}), we obtain
$$\nLdi{S(t)u_{0N}S(t)v_{0N}}{(0,1)\times B}^2 \leq c \min(N,L)^{\frac 23 + \frac 23 + \frac 12+\epsilon} \nLdi{u_{0N}}{B}^2 \nLdi{v_{0L}}{B}^2.$$ 
\end{proof}

\begin{rem}
\label{JustifGNL}
The bilinear estimate (\ref{estJnJp}) of Bessel functions $J_\nu(z_\nu r)$ is optimal, as it translates the qualitative behavior of $J_\nu$ near its first positif zero as an Airy function. The bilinear estimate on spherical harmonics \cite{BGTinv} is optimal, as can be checked on large degree spherical harmonics concentrating on large circles. Therefore, in order to obtain from (\ref{mainPsJn1}) a bilinear estimate $\Ps$ with $s<1$, we need to have a bound of $\Gamma_{NL}(\tau)$ by at most $\min(N,L)^{\gamma_0}$, for some $\gamma_0<\frac 56$.
\end{rem}

\begin{center}
    {\bf APPENDIX : proofs of two technical lemmas}
  \end{center}
\appendix

\section{Proof of Lemma \ref{cardGnl}}
\label{proofCardGnl}
The proof of this lemma uses the method of exponential sums and the precised asymptotic of $z_{n+\frac 12}$ from Lemma \ref{asymptzn1}. Hereafter we shall consider $N \leq L$. 

Let us remind the asymptotics of the first positif zero of Bessel function given by Lemma \ref{asymptzn1} : $z_\nu = \nu +\alpha \nu ^{\frac 13} + \beta \nu^{-\frac 13} + O(\nu^{-\frac 23}).$ Then $z_\nu^2 = \nu^2 +2\alpha \nu ^{\frac 43} + (2\beta +\alpha^2)\nu^{\frac 23} + O(\nu^\frac 13).$ For $t\in \RR_+$, let us denote by
\begin{equation}
\label{defg}
g(t) = t^2 +2\alpha t^{\frac 43} + (2\beta +\alpha^2)t^{\frac 23}.
\end{equation}
From the definition of $\Gamma_{NL}(\tau)$ (\ref{defGNL}), there exist $C_1,C_2>0$ such that $\Gamma_{NL}(\tau)$ is a subset of $\tilde{\Gamma}_{NL}(\tau)$, where $\tilde{\Gamma}_{NL}(\tau)$ is defined by
\begin{equation}
\label{slargset}
\{ (n,p)\in \mathbb{N}^2\ s.t. \ n\sim_\delta N,\ p\sim_\delta L,\ \tau -C_1 L^\frac 13 \leq g(n+\frac 12) + g(p+\frac 12) \leq \tau +C_2 L^\frac 13\}.
\end{equation}

We begin with a general result connecting $N$ and $L$. Let us recall that we have considered $N\leq L$.

\begin{lem}
\label{LN2}
For $N,L \in \delta^\mathbb{N}$ large, if $\Gamma_{NL}(\tau)$ contains at least two distinct elements, then there exists $c>0$ such that $L\leq c N^2$.
\end{lem}
\begin{proof}
Let us begin by considering $(n_1,p_1)\neq (n_2,p_2)$ two elements of $\Gamma_{NL}(\tau)$ such that $p_2 \neq p_1$. As $\Gamma_{NL}(\tau)\subset \tilde{\Gamma}_{NL}(\tau)$, we deduce from (\ref{slargset}) that $(n_1,p_1),(n_2,p_2)$ verify
$$|g(p_2+\frac 12) - g(p_1+\frac 12)| - cL^\frac 13 \leq  |g(n_2+\frac 12) - g(n_1+\frac 12)|.$$
From the definition of $g$ in (\ref{defg}), we have $g(p_2+\frac 12) - g(p_1+\frac 12) = (p_2-p_1) G(p_1,p_2)$, where $G(p_1,p_2) = p_1 + p_2+O(L^{\frac 13})$. As $|p_2-p_1|\geq 1$, we have $|G(p_1,p_2)| |p_2-p_1| \geq L$. On the other side, from (\ref{defg}) we deduce that $$|g(n_2+\frac 12) - g(n_1+\frac 12)| \leq cN^2.$$ Combining those estimates, we get $L \leq cN^2$. 

Let us explain how we obtain a contradiction in the case where all couples $(n_1,p_1)\neq (n_2,p_2) $ from $\Gamma_{NL}(\tau)$ are such that $p_1=p_2=p_0$. 
From (\ref{defGNL}), $\Gamma_{NL}(\tau)=\{ (n,p)\in \mathbb{N}^2\ s.t. \ n\sim_\delta N,\ \tau_0 - C_1 N^\frac 13 \leq g(n+\frac 12) \leq \tau_0 +C_2 N^\frac 13\},$ where we have denoted by $\tau_0=\tau-g(p_0 + \frac 12)$. We deduce $cN\leq |g(n_2+\frac 12) - g(n_1+\frac 12)| \leq cN^\frac 13$, which is absurd. 

Consequently, the case where all couples $(n_1,p_1)\neq (n_2,p_2) $ of $\Gamma_{NL}(\tau)$ are such that $p_1=p_2=p_0$ cannot occur. Thus, if $\Gamma_{NL}(\tau)$ contains two distinct elements, we can consider them such that $p_1\neq p_2$. As we showed previously this implies $L \leq cN^2$. 
\end{proof}

Let $x_{min}=g(N)+g(L)$ and $x_{max}=g(\delta N) + g(\delta L)$. We define a larger set than $\tilde{\Gamma}_{NL}(\tau)$ (\ref{slargset}) and consider its cardinal : for $x\in[x_{min}, x_{max}]$, $$S(x,N)=\#\{ (n,p)\in \mathbb{N}^2\ s.t. \ n\sim_\delta N,\ g(n+\frac 12) + g(p+\frac 12) \leq x \}.$$
For $\tau \in [x_{min}+C_1 L^\frac 13, x_{max}-C_2 L^\frac 13]$,
\begin{equation}
\label{SxNjustif}
\#\tilde{\Gamma}_{NL}(\tau) \leq S(\tau +C_2 L^\frac 13) - S(\tau -C_1L^\frac 13).
\end{equation}
Notice that for $x<x_{\min}$ the set that defines $S(x_{\min},N)$ minus the set that defines $S(x,N)$ contains only elements that are not in $\Gamma_{NL}(\tau)$, for $\tau=x+C_1 L^{\frac 13}$. Similarly for $x>x_{\max}$, the set that defines $S(x,N)$ minus the set that defines $S(x_{\max},N)$ contains only elements that are not in $\Gamma_{NL}(\tau)$, for $\tau=x-C_2 L^{\frac 13}$.

For $x\in [x_{min}, x_{max}]$ and $t\in [N,\delta  N]$ let us define $$f_x(t)=g^{-1}(x-g(t)) - \frac 12.$$
Therefore, 
\begin{equation}
\label{SxNf}
S(x,N)= \sum_{n\sim_\delta N} [f_x(n+\frac 12)] = \sum_{n\sim_\delta N} f_x(n+\frac 12) - \sum_{n\sim_\delta N} \{f_x(n+\frac 12)\},
\end{equation}
where we have denoted respectively by $[y]$ and $\{y\}$ the integer and the fractional part of $y$.

Using the definition of $g$ in (\ref{defg}), we obtain, via direct computation, informations on the behavior of $f_x(t)$.
\begin{lem} 
For $x\in [x_{min}, x_{max}]$ and $t\in [N,\delta N]$, 
\begin{equation}
\label{estFxt}
f_x(t) \sim L, \ f_x{'}(t) \sim -\frac NL, \ f_x{''}(t) \sim -\frac 1L, \ \partial_x f_x(t) \sim \frac 1L.
\end{equation}
\end{lem}
\begin{proof} Once we show that $f_x(t) \sim L$, the other estimates follow easily. From the definition of $g$ (\ref{defg}) we have that for large $N$ and $t\in [N,\delta N]$, $t\mapsto g(t)$ is a positif, increasing, convex function which behaves like $t\mapsto t^2$.  It suffices to find constants $c>0$ and $C>0$ such that $$cL^2 \leq x-g(t) \leq CL^2.$$ 

The right bound is easy. Since $x\leq x_{max}=g(\delta N) + g(\delta L)$, $$x-g(t) \leq g(\delta L) + g(\delta N) - g(N).$$ Using that for large $N$, $g(\delta N)\leq 1.1 \delta^2 N^2$, we obtain $x-g(t) \leq g(\delta L) + (1.1 \delta^2-1) N^2$ and the existence of $C>0$ such that $x-g(t)\leq CL^2$ follows easily from $N\leq L$. 

The proof of the inferior bound uses the fact that $N,L \in \delta^{\mathbb{N}}$ : $$x-g(t)\geq g(N)+g(L)-g(\delta N)\geq (1-1.1 \delta^2) N^2 + g(L)\geq (2-1.1 \delta^2) L^2.$$ Since $\delta>1$ and close to $1$, there exists $c>0$ such that $x-g(t)>cL^2$

\end{proof}

We now turn to the main part of the proof of Lemma \ref{cardGnl}, using a variant of the circle method, which can be found for example in \cite{Iwan}.

We begin with a general lemma on the fractional part. The fractional part is periodic of period $1$ and can be written as a Fourier series. The next lemma gives the exact convergence. Let us denote by $\nX{\theta}{}= \{\min |\theta - m|,\ m\in \mathbb{Z}\}.$

\begin{citelem}{\upshape (e.g.\cite{Iwan})} Let $H\geq 2$. For $\theta \in \RR$ one has
$$\{\theta\} = \frac 12 -\sum_{0<|h|< H} \frac{e^{2i\pi \theta h}}{2i\pi h} + O(\min(1, \frac 1{H\nX{\theta}{}})).$$
Moreover, $\min(1, \frac 1{H\nX{\theta}{}})$ is a periodic function of period $1$ and its Fourier series reads as follows :
$$\min(1, \frac 1{H\nX{\theta}{}}) = \sum_{h=-\infty}^{\infty} a_{h,H} e^{2i\pi \theta h}.$$
The coefficients $a_{h,H}$ satisfy the following estimates : there exist $c>0$ such that
$|a_{0,H}|\leq c \frac{ln H}{H}$, for  $0<|h|\leq H$ : $|a_{h,H}|\leq \frac c{H}$ and for $H\leq |h|$ : $|a_{h,H}|\leq c\frac{H}{|h|^2}.$
\end{citelem} 

Replacing these expansions into (\ref{SxNf}), the expression of $S(x,N)$ can be written as $S(x,N) = MT(x,N) + TE_1(x,N) + TE_2(x,N),$ where we have denoted by $$MT(x,N) = \sum_{n\sim_\delta N} (f_x(n+\frac 12) - \frac 12)$$ the main term of the sum and by $$TE_1(x,N) = \sum_{0<|h|<H} (\frac {1}{2i\pi h} + a_{h,H}) \sum_{n\sim_\delta N} e^{2i\pi f_x(n+\frac 12) h}$$ and by $$TE_2(x,N) = a_{0,H} N + \sum_{|h|\geq H} a_{h,H} \sum_{n\sim_\delta N} e^{2i\pi f_x(n+\frac 12) h}$$ the error terms. 
Let us recall a classical result, the Van der Corput Lemma, that allows us to estimate the trigonometric sums $\sum_{n\sim_\delta N} e^{2i\pi f_x(n+\frac 12) h}$ :

\begin{citelem}{\upshape (e.g.\cite{Iwan})} Let $b-a \geq 1$. Let $r(t)$ be a real function on $[a,b]$ such that $-\eta \Lambda \leq r{''}(t) \leq -\Lambda$ with $\Lambda>0$ and $\eta\geq 1$. Then there exists $c>0$ such that
$$|\sum_{a<n<b} e^{2i\pi r(n)}| \leq c (\eta \Lambda^{\frac 12}(b-a) + \Lambda^{-\frac 12}).$$
\end{citelem}

We want to apply this lemma with $a=N$, $b=\delta N$ and $r(t)=hf_x(t+\frac 12)$. From estimates (\ref{estFxt}) we deduce $\Lambda = c\frac h L$. Thus,
\begin{equation}
\label{vdCor}
|\sum_{n\sim_\delta N} e^{2i\pi f_x(n+\frac 12) h}| \leq c(Nh^{\frac 12}L^{-\frac 12} + h^{-\frac 12} L^{\frac 12}).
\end{equation}

First we treat $TE_2$. We have $|a_{0,H} N| \leq c NH^{-1}lnH.$ 
For $|h|\geq H$, we use $|a_{h,H}|\leq c \frac {H}{|h|^2}$ and (\ref{vdCor}) to conclude that
\begin{eqnarray*}
\left| \sum_{|h|\geq H} a_{h,H} \sum_{n\sim_\delta N} e^{2i\pi f_x(n+\frac 12) h} \right| &\leq & c\sum_{|h|\geq H}  H(Nh^{-\frac 32} L^{-\frac 12} + h^{-\frac 52}L^{\frac 12}) \\ 
&\leq & c (NH^{\frac 12} L^{-\frac 12} + H^{-\frac 12} L^{\frac 12}).
\end{eqnarray*}

We resume 
\begin{equation}
\label{estTE2}
|TE_2(x,N)|\leq c (NH^{-1}lnH + NH^{\frac 12} L^{-\frac 12} + H^{-\frac 12} L^{\frac 12}).
\end{equation}

If we apply the Van der Corput lemma to estimate $TE_1(x,N)$, we obtain $|TE_1(x,N)| \leq c (NH^\frac 12 L^{-\frac 12} + L^\frac 12)$. We know from Lemma \ref{LN2} that $L\leq cN^2$. Therefore this gives a bound of at least $cN$ for $\#\Gamma_{NL}(\tau)$, which is larger than the one we announced in Lemma \ref{cardGnl}. We shall proceed differently. 

As we want to estimate $S(x_2,N)-S(x_1,N)$, we shall evaluate the difference for each term of the decomposition $$S(x,N) = MT(x,N) + TE_1(x,N) + TE_2(x,N). $$ For the main term this reads
$$MT(x_2,N)-MT(x_1,N) = \sum_{n\sim_\delta N} f_{x_2}(n+\frac 12) - f_{x_1}(n+\frac 12).$$
From the mean value theorem and estimates (\ref{estFxt}) we deduce $$|MT(x_2,N)-MT(x_1,N)| \leq c N(x_2-x_1)L^{-1}.$$
We pass to the estimate of $TE_1(x_2,N)-TE_1(x_1,N)$ : $$\sum_{0<h<H,\ n\sim_\delta N}(\frac {1}{h} + ca_{h,H}) \left( \sin(2\pi h f_{x_2}(n+\frac 12))- \sin(2\pi h f_{x_1}(n+\frac 12)) \right).$$ From the mean value theorem and estimates (\ref{estFxt}) we have that the difference of sines above is bounded by $ch |x_2-x_1| \sup_{t\in [N,2N]} \partial_x{f_x}(t) \leq ch L^{-1} |x_2-x_1|.$ We multiply it by $|\frac 1h + ca_{h,H}| \leq \frac{c}{h}$ and sum over $0<h<H$ and $n\sim_\delta N$. Consequently,
\begin{equation}
\label{VIP}
|TE_1(x_2,N)-TE_1(x_1,N) | \leq c NHL^{-1}|x_2-x_1|.
\end{equation}
Combining these with the bounds on $TE_2$ (\ref{estTE2}) we conclude that $|S(x_2,N)-S(x_1,N)| \leq F(H)$, where
\begin{equation}
\label{VVIP}
 F(H)=NHL^{-1}|x_2 -x_1| + NH^{-1+\epsilon} + NH^{\frac 12} L^{-\frac 12} + H^{-\frac 12} L^{\frac 12}.
\end{equation}

From (\ref{SxNjustif}) we know that we have to consider $x_2 = \tau +C_2 L^\frac 13$ and $x_1=\tau - C_1 L^\frac 13$. Consequently, $x_2-x_1=CL^\frac 13$ and thus $F(H)=NHL^{-\frac 23} + NH^{-1+\epsilon} + NH^{\frac 12} L^{-\frac 12} + H^{-\frac 12} L^{\frac 12}$. We choose $H$ that minimizes $F(H)$ : $H=L^{\frac 13}$ satisfies and gives $F(L^\frac 13)=cNL^{-\frac 13 + \epsilon} + L^{\frac 13}$. Using the hypothesis $N\leq L$ and $L<cN^2$ (see Lemma \ref{LN2}), this choice of $H$ yields the announced result
$$\#\Gamma_{NL}(\tau)\leq S(\tau + C_2 L^\frac 13,N)-S(\tau + C_1 L^\frac 13,N) \leq cN^{\frac 23+\epsilon}.$$

\section{Proof of Lemma \ref{asymptzn1}}
\label{pfAsymptzn1}

Let us begin by motivating the need of a precise asymptotics of $z_\nu$, the first positif zero of the Bessel function $J_\nu$. As we pointed out in Remark \ref{JustifGNL}, we need to prove a bound on $\Gamma_{NL}(\tau)$ defined by (\ref{defGNL}) of at most $\min(N,L)^{\gamma_0}$ for some $\gamma_0<\frac 56$.
In (\ref{SxNjustif}) we bounded $\#\Gamma_{NL}(\tau) \leq S(\tau+C_2 L^\frac 13) - S(\tau+C_1 L^\frac 13)$. If we want to make a similar analysis to the one done in the proof of Lemma \ref{cardGnl} (Appendix \ref{proofCardGnl}) for the Sturm asymptotics $z_\nu = \nu + \alpha \nu^{\frac 13} + O(\nu^{-\frac 13})$, we have to consider the equivalent of $g$ defined in (\ref{defg}) : $\bar{g}(t) = t^2 + 2\alpha t^\frac 43$ and the equivalent of $\tilde{\Gamma}_{NL}(\tau)$ from (\ref{slargset}) : $\bar{\Gamma}_{NL}(\tau)=\{ (n,p)\in \mathbb{N}^2\ s.t. \ n\sim N,\ p\sim L,\ \tau +C_1 L^\frac 23 \leq \bar{g}(n+\frac 12) + \bar{g} (p+\frac 12) \leq \tau +C_2 L^\frac 23\}.$ Notice that the estimates on $\bar{g}$ are not fundamentally influenced by terms of order less than $2$. Same analysis as the one on $g$ and $\Gamma_{NL}(\tau)$ can be performed. Thus, (\ref{VVIP}) becomes in these setting 
$$\bar{\Gamma}_{NL}(\tau) \leq c(NHL^{-\frac 13} + NH^{-1+\epsilon} + NH^{\frac 12} + H^{-\frac 12} L^\frac 12).$$ Using $N\leq L$ and $L\leq cN^2$, for every choice of $H$ we obtain a bound larger than the $\min(N,L)^{\frac 56}$ needed in order to obtain a bilinear estimate with $s<1$ (see Remark \ref{JustifGNL}). Consequently, having precise asymptotics of $z_\nu$ helps us to give precise bounds on the distribution of all first zeros of Bessel functions.

In order to simplify the proof, we make use of the Sturm asymptotics $z_\nu = \nu + \alpha \nu^{\frac 13} + O(\nu^{-\frac 13})$. We follow a similar approach to the one used in \cite{BGTcortona}. That is, we give a precise asymptotics of $J_\nu$ close to its first zero in terms of Airy function. We identify highest order terms and invert the asymptotics. We inspire the proof of the asymptotics from that of Theorem 7.7.8 from \cite{Ho}. The mentioned theorem gives abstract asymptotics of an oscillatory integral with a phase with critical point of order $2$. This is exactly the case of the Bessel function near its first positif zero and moreover we can compute the first terms of the expansion.

We now pass to the proof itself. We have defined $J_\nu(x)=T_{1,\nu}(x) + T_{2,\nu}(x)$ in (\ref{defJnu}), where $T_{1,\nu}(\nu \rho)=\frac{1}{2\pi} \int_{-\pi}^\pi e^{i\nu \Phi_\rho(\theta)} d\theta$ and $\Phi_\rho(\theta)=\theta - \rho \sin \theta$. We write it under the form $T_{1,\nu}(\nu \rho)=\frac{1}{2\pi} \int_{-\pi}^\pi e^{i\nu ((1-\rho)\theta + \rho \frac{\theta^3} {6})} e^{i\nu \rho r_1(\theta)} d\theta$, where we denote by $r_1(\theta) = \sin \theta - \theta + \frac {\theta^3}{6}.$ For $\theta$ in the neighborhood of $0$ we have $r_1(\theta)= O(\theta^5)$. We perform the change of variable $\rho \frac{\theta^3} {2} = \xi^3$. We denote by $$a(\rho) = (1-\rho) \left( \frac{2}{\rho}\right)^\frac 13$$ and by $\tilde{r}_1(\xi) = r_1((\frac{2}{\rho})^\frac 13 \xi)$. Thus,
$$T_{1,\nu}(\nu \rho)=\frac{1}{2\pi} \left( \frac{2}{\rho}\right)^\frac 13 \int_{-\pi \left( \frac{\rho}{2}\right)^\frac 13} ^{\pi \left( \frac{\rho}{2} \right)^\frac 13 } e^{i\nu (a(\rho) \xi + \frac{\xi^3} {3})} v_\nu(\rho,\xi) d\theta,$$
where $v_\nu(\rho, \xi) = e^{i\nu \tilde{r}_1(\xi)}.$ 
Let us recall that using the Sturm asymptotics of $z_\nu=\nu \rho_0$, we know that $\rho_0 = 1 + \alpha \nu^{-\frac 13} + O(\nu^{-\frac 43})$. Therefore, we are interested in the behavior of $T_{1,\nu}(\nu \rho)$ for $\rho = 1 +O(\nu^{-\frac 13})$.
Let us consider $\chi \in C_0^\infty(\RR)$ such that $\chi(x)=1$ for $|x|\leq 1$ and $\chi(x)=0$ for $|x|\geq 2$. Let $\lambda>0$ to be chosen such that $T_{1,\nu}(\nu \rho)=\frac{1}{2\pi} \left( \frac{2}{\rho}\right)^\frac 13 \int_{\RR} e^{i\nu (a(\rho) \xi + \frac{\xi^3} {3})} \chi(\nu^{\lambda} \xi) v_\nu(\rho,\xi) d\theta + TR_\nu$, where $TR_\nu$ small. Applying the nonstationary phase estimate we obtain that for  $|\rho-1|< c \nu^{-\frac 23}$ and $\lambda<\frac 14$, $|TR_\nu|\leq c_k \nu^{-k}$. 

By the division theorem (Weierstrass formula) (e.g. Theorem 7.5.2, \cite{Ho}), we have the existence and uniqueness of $q_\nu(\rho,\xi)$, $R_{1,\nu}(\rho)$ and $R_{2,\nu}(\rho)$ such that
$$\chi(\nu^\lambda \xi) v_\nu(\rho, \xi) = (\xi^2 + a(\rho)) q_\nu(\rho,\xi) + R_{2,\nu}(\rho) \xi  + R_{1,\nu}(\rho).$$ 
Let $\xi_+ = \sqrt{-a(\rho)}$. Then $R_{1,\nu}(\rho) = \cos(\nu \rho \tilde{r}_1(\xi_+)) \chi(\nu^\lambda \xi_+)$ and $R_{2,\nu}(\rho) = i \frac{\sin(\nu \rho \tilde{r}_1(\xi_+))}{\xi_+} \chi(\nu^\lambda \xi_+)$. Using $|\rho-1|\leq c\nu^{-\frac 23}$, $a(\rho)=(1-\rho)\left( \frac{2}{\rho} \right)^{\frac 13}$ and $\lambda<\frac 14$, we obtain $\chi(\nu^\lambda \xi_+)=1$ and consequently 
\begin{equation}
\label{R1R2}
R_{1,\nu}(\rho) = \cos(\nu \rho \tilde{r}_1(\xi_+)) \ {\rm and} \ R_{2,\nu}(\rho) = i \frac{\sin(\nu \rho \tilde{r}_1(\xi_+))}{\xi_+}.
\end{equation}

Let us notice that $$\int e^{i\nu \left( a(\rho)\xi + \frac{\xi^3}{3} \right)}d\xi = \nu^{-\frac 13} Ai(\nu^{\frac 23} a(\rho))$$ and $$\int e^{i\nu \left( a(\rho)\xi + \frac{\xi^3}{3} \right)} \xi d\xi = - i \nu^{-\frac 23} Ai'(\nu^{\frac 23} a(\rho)).$$
Thus,
\begin{eqnarray*}
T_{1,\nu}(\nu \rho)= \left( \frac 2 \rho \right)^{\frac 13} \int e^{i\nu (a(\rho)\xi + \frac {\xi^3} {3})} \left( a(\rho)+ \xi^2 \right) q_\nu(\rho, \xi) d\xi \\
+\left( \frac 2 \rho \right)^{\frac 13} i^{-1} R_{2,\nu}(\rho) \nu ^{-\frac 23} Ai'(\nu^{\frac 23} a(\rho)) \\
+\left( \frac 2 \rho \right)^{\frac 13} R_{1,\nu}(\rho) \nu ^{-\frac 13} Ai(\nu^{\frac 23} a(\rho)) 
\end{eqnarray*}

Let us denote the first term of the sum by $I_{1,\nu}$. In this term we can perform an integration by parts and we obtain
$$I_{1,\nu} = \left( \frac 2 \rho \right)^{\frac 13} \nu^{-1} \int e^{i \nu (a(\rho)\xi + \frac {\xi^3} {3})} \partial_\xi q_{\nu}(\rho, \xi) d\xi.$$
 Using the integral representation of $q_\nu(\rho, \xi)$ from (7.5.4) \cite{Ho}, we have, for $|\rho -1|< c\nu^{-\frac 23}$ and $\xi<c$,
$$q_\nu(\rho, \xi) = \frac{1}{2\pi i} \int_{|\eta-1|=c\nu^{-\frac 23}} \frac{\chi(\nu^\lambda \xi) v_\nu(\eta, \xi)} {(\xi^2 + a(\eta))(\eta-\rho)} d\eta.$$ We obtain that $|I_{1,\nu}| \leq c\nu^{-1}$ (in fact we have better estimates, but we only use this one since $|T_{2,\nu}| \leq c \nu^{-1}$).

Putting together the estimates on $T_{1,\nu}(\nu \rho)$ and $T_{2,\nu}(\nu \rho)$, we obtain that $J_{\nu}(\nu \rho)$ equals
\begin{equation}
\label{Jnuz0}
\left( \frac 2 {\nu \rho} \right)^{\frac 13} R_{1,\nu}(\rho)  Ai(\nu^{\frac 23} a(\rho)) + \left( \frac 2 {\nu^2 \rho} \right)^{\frac 13} i^{-1} R_{2,\nu}(\rho) \nu ^{-\frac 13} Ai'(\nu^{\frac 23} a(\rho)) +O(\nu^{-1}).
\end{equation}

Let us recall that we want to give a better asymptotics for $z_\nu=\nu \rho_0$ of the first positif zero of Bessel function $J_\nu$. From Sturm asymptotics (\cite{Wat}, chap. XV) we have that $\rho=1-2^{-\frac 13} z_0 \nu^{-\frac 23} + O(\nu^{-\frac 43})$, where $z_0$ is the first negatif zero of the Airy function. Thus, $a(\rho)\nu^{\frac 23} - z_0 = O(\nu^{-\frac 23})$. Using Taylor expansion of the Airy function near its first negatif zero and the identity verified by the Airy function $Ai''(x)-xAi(x)=0$, we have $$Ai(\nu^{\frac 23} a(\rho)) = (a(\rho)\nu^{\frac 23} - z_0)Ai'(z_0) + O((a(\rho)\nu^{\frac 23}- z_0)^3)$$ and $$Ai'(\nu^{\frac 23}a(\rho)) = Ai'(z_0) + O((a(\rho)\nu^{\frac 23} - z_0)^2).$$ From (\ref{R1R2}) we deduce that, for $\rho=1+O(\nu^{-\frac 23})$, we have $\nu \rho \tilde{r}_1(\xi_+)=O(\nu^{-\frac 23})$ and therefore $R_{1,\nu}(\rho) = 1 + O(\nu^{-\frac 43})$. Moreover, $i^{-1} R_{2,\nu}(\rho)= \nu \rho \left( \frac 2 \rho \right)^{\frac 53} \frac {a(\rho)^2}{5!} + O(\nu^{-1})$. 
We replace those expansions in the equation $J_{\nu}(\nu \rho_0) = 0$ using the expression of $J_{\nu}(\nu \rho)$ from (\ref{Jnuz0}). We identify the highest order terms and obtain
$$a(\rho) \nu^{\frac 23} - z_0 + \frac{2^\frac 53}{5!} z_0^2 \nu^{-\frac 23} = O(\nu^{-1}).$$
Consequently, $a(\rho) = z_0 \nu^{-\frac 23} + \frac{2^\frac 53}{5!} z_0^2 \nu^{-\frac 43} + O(\nu^{-\frac 5 3}).$ We write $\rho=1+\gamma$, $|\gamma|\leq c \nu^{-\frac 23}$. We have $a(1+\gamma) = -\gamma \left( \frac{2}{1+\gamma}\right)^\frac 13$. Let us denote this function by $b(\gamma)$. We compute $b^{-1}(r)$ for $r$ in a neighborhood of zero. Inverting the asymptotics we obtain $b^{-1}(r)=-2^{-\frac 13} r -\frac 13 r^2 + O(r^3)$. As $\gamma= b^{-1}(z_0 \nu^{-\frac 23} + \nu^{-\frac 43} \frac{2^\frac 53}{5!} z_0^2 + O(\nu^{-\frac 5 3}))$ we obtain $\gamma = -2^{-\frac 13} \nu^{-\frac 23} z_0 + c z_0^2 \nu^{-\frac 43} + O(\nu^{-\frac 53})$. This yields $z_\nu=\nu(1+\gamma)=\nu -2^{-\frac 13}z_0 \nu^\frac 13 + c z_0^2 \nu^{-\frac 13} +O(\nu^{-\frac 23})$.


\end{document}